\documentclass[a4paper,fleqn]{cas-sc}

\usepackage[numbers, sort&compress]{natbib}

\usepackage{amsmath,amsopn,amsthm,amssymb}
\usepackage{graphicx}
\usepackage[mathscr]{euscript}
\usepackage{mathtools}
\usepackage{float}
\usepackage{enumitem}
\usepackage{xspace}
\usepackage{scalerel}
\usepackage{marginnote}
\usepackage{enumitem}

\usepackage{algorithmicx}
\usepackage{algorithm} 
\usepackage[noend]{algpseudocode}
\algrenewcommand\algorithmicrequire{\textbf{Input:}}
\algrenewcommand\algorithmicensure{\textbf{Output:}}

\newtheorem{theorem}{Theorem}
\newtheorem{lemma}[theorem]{Lemma}

\newtheorem{proposition}[theorem]{Proposition}
\newtheorem{corollary}[theorem]{Corollary}
\newtheorem{example}[theorem]{Example}
\newtheorem{definition}[theorem]{Definition}
\newtheorem{observation}[theorem]{Observation}
\newtheorem*{remark}{Remark}

\newcommand{\mh}[1]{\begingroup#1\endgroup}
\newcommand{\REV}[1]{\begingroup#1\endgroup}
\newcommand{\REVII}[1]{\begingroup#1\endgroup}

\newcommand{\cl}{\ensuremath{\operatorname{cl}}}
\newcommand{\tc}{\ensuremath{\operatorname{tc}}}
\newcommand{\lca}{\ensuremath{\operatorname{lca}}}
\newcommand{\LCA}{\ensuremath{\operatorname{LCA}}}
\newcommand{\Hasse}{\mathscr{H}} 
\newcommand{\notR}{\mathrel{R\!\!\!/}}
\newcommand{\cG}{\mathcal{G}} 
\newcommand{\cN}{\mathcal{N}} 

\newcommand{\srel}{\smallblacktriangleleft}

\newcommand{\rel}{\trianglelefteq}
\newcommand{\altsrel}{S}

\newcommand{\axiom}[1]{\textnormal{\textbf{(#1)}}}
\newcommand{\pairs}{\mathcal{P}_2}
\newcommand{\support}{\ensuremath{\operatorname{supp}}}

\DeclareMathOperator{\indeg}{indeg}
\DeclareMathOperator{\outdeg}{outdeg}
\DeclareMathOperator{\CC}{\mathtt{C}}

    {\begin{description}[leftmargin = 0.2cm, labelsep = 0.2cm]}
    {\end{description}}

\def\tsc#1{\csdef{#1}{\textsc{\lowercase{#1}}\xspace}}
\tsc{AL}
\tsc{AA}
\tsc{VM}
\tsc{GES}
\tsc{MH}

\begin{document}
\let\WriteBookmarks\relax
\def\floatpagepagefraction{1}
\def\textpagefraction{.001}

\shorttitle{Inferring DAGs and Phylogenetic Networks from Least Common Ancestors}  
\shortauthors{{Lindeberg et~al.}}  

\title[mode = title]{Inferring DAGs and Phylogenetic Networks from Least Common Ancestors}

\author[1]{Anna Lindeberg}[orcid=0000-0001-9664-1918]
\ead{anna.lindeberg@math.su.se}
\credit{conceptualized the main ideas. Developed and proved the theoretical concepts and wrote the main manuscript. Reviewed, edited, and approved the final manuscript}
	
\author[1]{Anton Alfonsson}[]
\ead{anton@alfonsson.se}
\credit{developed and proved the theoretical concepts. Implemented Algorithm~\ref{alg:REAL}. Reviewed, edited, and approved the final manuscript}

\author[2]{Vincent Moulton}[orcid=0000-0001-9371-6435]
\ead{v.moulton@uea.ac.uk}
\credit{developed and proved the theoretical concepts and wrote the main manuscript. Reviewed, edited, and approved the final manuscript}

\author[3]{Guillaume E. Scholz}[orcid=0000-0001-5033-8040]
\ead{gllm.scholz@gmail.com}
\credit{developed and proved the theoretical concepts and wrote the main manuscript. Reviewed, edited, and approved the final manuscript}

\author[1]{Marc Hellmuth}[orcid=0000-0002-1620-5508]
\cormark[1]
\ead{marc.hellmuth@math.su.se}
\credit{conceptualized the main ideas. Developed and proved the theoretical concepts and wrote the main manuscript. Reviewed, edited, and approved the final manuscript}

\cortext[1]{Corresponding author}

\affiliation[1]{organization={Department of Mathematics, Stockholm University},
            addressline={Albanov{\"a}gen 28}, 
            city={Stockholm},
            postcode={SE-106 91}, 
            country={Sweden}}

\affiliation[2]{organization={School of Computing Sciences, University of East Anglia},
            addressline={Norwich Research Park}, 
            city={Norwich},
            postcode={NR4 7TJ}, 
            country={United Kingdom}}

\affiliation[3]{organization={Independent Researcher},
            city={Leipzig},
            postcode={DE-04229}, 
            country={Germany}}


\begin{abstract}
A least common ancestor (LCA) of two leaves in a directed acyclic graph (DAG) is a vertex  
that is an ancestor of both of the leaves, and for which no proper descendant 
of this vertex is also an ancestor of the two leaves. 
LCAs play a central role in representing hierarchical relationships in rooted trees and, more generally, in DAGs. 
In 1981, Aho et al.\ introduced the problem of determining whether a collection of pairwise LCA constraints
on a set $X$, of the form $(i,j)<(k,l)$ with $i,j,k,l \in X$,
can be realized by a rooted tree with leaf set $X$ 
such that, whenever $(i,j) < (k,l)$ holds, the LCA of $i$ and $j$ in the tree is a descendant of the LCA of $k$ and $l$. 
In particular, they presented a polynomial-time algorithm, called \textsc{Build}, to solve this problem.
In many cases, however, such constraints cannot be realized by \emph{any} tree. 
In these situations, it is natural to ask whether they can still be realized by a more general directed acyclic graph (DAG).  

In this paper, we extend Aho. et~al.'s problem from trees to DAGs, providing a 
theoretical and algorithmic framework for reasoning about LCA constraints in this more general setting. 
More specifically, given a collection $R$ of LCA constraints, we introduce the notion of the 
\emph{$+$-closure} $R^+$, which captures additional LCA relations implied by $R$. 
Using this closure, we then associate a \emph{canonical DAG} $\cG_R$ to $R$ and show 
that a collection of constraints $R$ can be realized by a DAG in terms of LCAs if and only if $R$ is realized in this way by $\cG_R$. 
In addition, we adapt this construction to obtain {\em phylogenetic networks}, a special class of DAGs widely used in phylogenetics. 
In particular, we define a canonical phylogenetic network $\cN_R$ for realizing $R$, and prove that $\cN_R$ is \emph{regular}, 
that is, it coincides with the Hasse diagram of the set system underlying $\cN_R$.
Finally, we show that, for any DAG-realizable collection $R$, the 
{\em classical closure} of $R$, that is the collection of \emph{all} LCA constraints that must hold in \emph{every} DAG realizing $R$,
coincides with its $+$-closure. 
All constructions introduced in this paper can be computed in polynomial time, and we provide explicit algorithms for each.
\REV{All  algorithms developed in this paper are implemented in the freely available Python package \texttt{RealLCA}.}
\end{abstract}





\maketitle

\sloppy
\section{Introduction}

\begin{figure}
    \centering
    \includegraphics[width=0.85\textwidth]{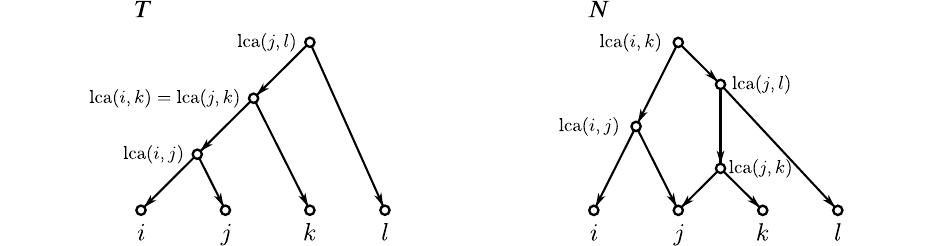}
    \caption{
        A tree $T$ and a network $N$, with least common ancestors indicated next to each vertex. The LCA constraints $(i,j)<(i,k)$ and $(j,k)<(j,l)$ are realized by both $T$ and $N$. The constraints $(i,j)<(i,k)$, $(j,k)<(j,l)$ and $(j,l)<(i,k)$ cannot be realized by any tree, but are realized by $N$. 
        }
    \label{fig:intro}
\end{figure}

In a rooted tree, the least common ancestor of two vertices $x$ and $y$, 
denoted $\lca(x,y)$, is the vertex that is an ancestor of 
both $x$ and $y$, such that no proper descendant of this vertex is also an ancestor of both $x$ and $y$. In the seminal paper \cite{Aho:81}, Aho et al.\ studied the  
following problem that is related to phylogenetic trees and 
relational databases. Suppose we are given a collection of constraints of the form $(i,j) < (k,l)$, with $i \neq j$, $k \neq l$, and $i,j,k,l \in X$.  Can we construct a rooted tree $T$ with leaf set $X$ such that, for every constraint $(i,j) < (k,l)$,  $\lca(i,j)$ is a descendant of $\lca(k,l)$ in $T$, or determine that no such tree exists?
For example, the tree $T$ as shown in Figure~\ref{fig:intro} realizes the constraints
$(i,j)<(i,k)$ and $(j,k)<(j,l)$ since the vertex $\lca(i,j)$ is a proper descendant of $\lca(i,k)$, and 
$\lca(j,k)$ is a proper descendant of $\lca(j,l)$.

In \cite{Aho:81} Aho et al.\ presented an efficient algorithm for solving this problem, called
\textsc{Build}. Since the \textsc{Build} algorithm was introduced it has become an indispensable tool within phylogenetic studies, where \textsc{Build} and its variants 
have been used to construct trees representing the evolutionary history of genes, species, or other taxa. Several practical implementations and improvements have been developed over the years~\cite{Jansson:05,DF:16,Henzinger:99,Holm:01}.
In particular, constraints of the form $(i,j) < (k,l)$ can often be inferred directly from genomic sequence data, making such LCA-based constraints a useful approach to reconstructing phylogenetic trees~\cite{Stadler2020,Geiss2019,Geiss2020}.

In many cases, however, such constraints cannot be realized by any tree. 
For example, consider the constraints $(i,j) < (i,k)$, $(j,k) < (j,l)$, and $(j,l) < (i,k)$. 
It can be shown (e.g., using \textsc{Build}) that no tree can satisfy all three constraints simultaneously. 
Nevertheless, these constraints can be realized by the directed acyclic graph (DAG) $N$ illustrated in Figure~\ref{fig:intro}.
Surprisingly, and to the best of our knowledge, no prior work has studied
LCA constraints in the broader setting of DAGs.
A major challenge is that, in DAGs, the LCA of two leaves is not necessarily unique, and in some cases it may not exist at all
(see Figure~\ref{fig:exmpl-LCA} below for an illustrative example). 
Nonetheless, in circumstances where these issues do not arise, it becomes an interesting problem to determine whether a DAG  can realize a given collection of LCA constraints, for instance, in situations where they cannot be realized by a tree.

In this paper, we characterize LCA constraints that can be realized by DAGs and provide 
polynomial-time algorithms to determine their realizability and, when this is possible, to construct DAGs that realize them. 
In doing so, we also address several interesting combinatorial problems that arise when defining the closure of a collection of LCA constraints, that is, the collection obtained by repeatedly generating new constraints implied by the existing ones.

We now briefly outline the contents of the rest of this paper. After 
presenting the main definitions in the next section, 
Section~\ref{sec:realizabilty} formalizes what it means for a DAG to realize a set of LCA constraints, where we introduce two main variants of realizability.  
The first, \emph{strict realization}, generalizes the approach of Aho et al.~\cite{Aho:81}. 
Specifically, for $i,j,k,l \in X$, it assumes that the vertices $\lca(i,j)$ and $\lca(k,l)$ are both well-defined in the DAG $G$ realizing the given collection of constraints, and that if $(i,j)<(k,l)$, then $\lca(i,j)$ is a descendant of $\lca(k,l)$ in $G$.  
The second, more general notion, \emph{realization}, allows LCA pairs to coincide when constraints are ``partially symmetric'', while still respecting the partial order implied by asymmetric constraints.

We then address the problem of characterizing when a collection $R$ of LCA constraints is realizable by some DAG. To this end, in Section~\ref{sec:closure} we introduce the \emph{$+$-closure} operator for a collection $R$, an operator that produces a new collection $R^+$ that satisfies three simple rules, the most important being that $R^+$ is cross-consistent.  
Intuitively, cross-consistency captures the idea that certain ``local'' LCA relationships imply others. Specifically,
it states that if in a DAG $G$ we have that 
$\lca(i,k)$ and $\lca(j,l)$ are descendants of $\lca(x,y)$, then, whenever $\lca(i,j)$ is well-defined in $G$, it must also be a descendant of $\lca(x,y)$.
The $+$-closure $R^+$ therefore captures additional information about any DAG $G$ that could potentially realize $R$.  In Section~\ref{sec:canonDAG}, using the $+$-closure, we define an equivalence relation $\sim_{R^+}$ on pairs of elements in $X$ and define the \emph{canonical DAG} $\cG_R$ associated to $R$, whose vertices are the equivalence classes of $\sim_{R^+}$.  
In Theorem~\ref{thm:char}, we then show that $R$ is realizable by a DAG if and only if it is realized by $\cG_R$, and also show that this is equivalent to $R$ satisfying two simple axioms that are expressed in terms of $R^+$ (see Definition~\ref{def:X12}).

In Section~\ref{sec:canon-N}, we study the problem of realizing collections of LCA constraints using \emph{phylogenetic networks}, or simply \emph{networks}.  
These are a special class of DAGs that generalize phylogenetic trees and are commonly used to represent reticulate evolution (see, e.g., \cite[Chapter 10]{steel2016phylogeny}).  
Several important classes of networks have the property that every internal vertex with at least two children is the LCA of at least two leaves. This includes, for example, normal networks \cite{Willson2010} (a leading class of networks \cite{Francis:25}), as well as regular and level-1 networks (see, e.g., \cite{Hellmuth2023,HL:24,HL:26} for more details on LCAs in networks).  By 
modifying the canonical DAG $\cG_R$ of a realizable collection of constraints $R$, we define in Section~\ref{sec:canon-N} a canonical \emph{network} $\cN_R$ that realizes $R$.  
We show that $\cN_R$ is \emph{regular}, meaning that is closely related to the Hasse diagram of the sets formed by taking the descendants in $X$ of each vertex in $\cN_R$.  
Exploiting the fact that $R^+$ can be computed in polynomial time by repeatedly applying a set of implication rules (Theorem~\ref{thm:rules} \REV{and \ref{thm:R+-polytime}}), we also present a polynomial-time algorithm  for computing $\cN_R$ (see Algorithm~\ref{alg:REAL}).

A natural way to define a closure for realizable relations is to follow the approach used for triplets or quartets in phylogenetic trees \cite{SH:18,MaayanLevy24,Bryant97,BS:95}.  
In particular, given a realizable relation $R$, we define its ``classical'' closure, denoted by $\cl(R)$, to be the collection of LCA constraints that must hold in \emph{every} DAG realizing $R$.  
In Section~\ref{sec:proof-closure}, we show the somewhat surprising result that, for any realizable $R$, the classical closure coincides with the $+$-closure; that is, $\cl(R) = R^+$.  
In Section~\ref{sec:incomp}, we extend the realizability problem by seeking a DAG that satisfies a collection of constraints $R$ while simultaneously satisfying another set $S$ of incomparability constraints.  
Finally, in Section~\ref{sec:outlook}, we discuss 
several interesting possible directions for future research.

All  algorithms developed in this paper are implemented in the freely available 
Python package \texttt{RealLCA} \cite{github-AL}.


\section{Basics}
\label{sec:basic}

\paragraph{Sets and Relations.}
In what follows, $X$ will always denote a finite non-empty set. 
We denote with $\mathcal{P}(X)$ the powerset of $X$. A \emph{set system} $\mathfrak{C}$ on $X$ is a subset of 
 $\mathcal{P}(X)$.  
A set system $\mathfrak{C}$ on $X$ is \emph{grounded} if $\{x\}\in \mathfrak{C} $ for all $x\in X$ and
$\emptyset\notin \mathfrak{C}$, while $\mathfrak{C}$ is a \emph{clustering system} if it is grounded
and satisfies  $X\in \mathfrak{C}$. 
We let  $\pairs(X)\coloneqq\{\{a,b\} \mid a,b\in X\}$ (with $a=b$ allowed) denote the  set 
system consisting of all 1- and 2-element subsets of $X$.
We will often write $ab$, respectively, $aa$ for elements $\{a,b\}$, respectively, $\{a\}$ in $\pairs(X)$. Thus, $ab=ba$ always holds. 

Given a set $A$, a  subset $R\subseteq A\times A$ is a \emph{binary relation (on $A$).}
We shall often write $a\ R\ b$ instead of $(a,b)\in R$ and  $a\ \notR\ b$ instead of $(a,b)\notin R$, for $a,b \in A$. In addition, we sometimes write $p\ R\ \dots\ R\ q$ if there is a \emph{$(p,q)$-chain in $R$}, i.e., some $a_0,\dots,a_k$, $k\geq 1$ such that $a_0 = p\ R\ a_1\ R\  \dots\ R\ a_{k-1} R\ a_k=q$.

\begin{remark}
As all relations considered in this work are binary, we shall simply refer to them as \emph{relations}.
\end{remark}

Furthermore, we define the \emph{support $\support_R$} of a   relation $R$ on $A$ as
\[\support_R \coloneqq \{p\in A \mid \text{ there is some } q\in A \text{ with }  p\ R\ q \text{ or } q\ R\ p \},  \]
that is, the subset of $A$ that contains precisely those $p\in A$ that are in $R$-relation with some $q\in A$. We often consider relations $R$ on $A =\pairs(X)$ in which case we extend
$\support_R$ to obtain
\[\support_{R}^+ \coloneqq \support_{R}\ \cup\ \{xx\mid x\in X\}.\]

\begin{example}\label{exmpl1}
The elements in the   relation $R = \{(\{x\},\{a,b\}), (\{y\},\{a,b\}), (\{u,v\},\{a,b\})\} \subseteq \pairs(X) \times \pairs(X)$ 
with $X = \{x,y,a,b,u,v,w\}$ can alternatively be written as  $xx \ R\ ab$, $yy\ R\ ab$ and $uv\ R\ ab$. 
Moreover, $\support_R = \{xx, yy, uv,ab\}$ and $\support_{R}^+  = \{ xx, yy, aa,bb, uu,vv, ww, uv,ab \}$.
\end{example}

Let $R$ be a relation on $A$.
Then, $R$ is \emph{asymmetric} if $p\ R\ q$ implies $q\ \notR\ p$ for all $p,q\in A$,
and it is \emph{anti-symmetric} if $p\ R\ q$ and $q\ R\ p$ implies $p=q$ for all $p,q\in A$.
Moreover, $R$ is \emph{transitive}, if  $p\ R\ q$ and  $q\ R\ r$ implies $p\ R\ r$
for all $p,q,r\in A$. 
If $B\subseteq A$, we say that $R$ is \emph{$B$-reflexive} if $(b,b)\in R$ for every $b\in B$. 
A relation  $R$  on  $A$ is \emph{reflexive} if it is $A$-reflexive. \REV{
A \emph{poset} $(Q, \leq)$ is a set $Q$  equipped with a 
partial order $\leq$, i.e., a relation $\leq$ on $Q$
that is reflexive, transitive, and anti-symmetric.}

A key concept in this paper is that of closure operators, see e.g.\ \cite{CASPARD2003241,BS:95,SH:18}.
A \emph{closure operator} on a set $S$ is a map $\phi\colon \mathcal{P}(S) \to\mathcal{P}(S) $ that satisfies the following three properties for all $R,R'\in \mathcal{P}(S)$:
\begin{itemize}[noitemsep]
    \item[] \emph{Extensivity}: $R \subseteq \phi(R)$.
    \item[] \emph{Monotonicity}: $R \subseteq R'$ implies $\phi(R) \subseteq \phi(R')$,
    \item[] \emph{Idempotency}: $\phi(\phi(R)) = \phi(R)$. 
\end{itemize}

We let $\tc(R)$ denote the \emph{transitive closure} of a relation $R$, 
that is, the relation where $(p,q) \in \tc(R)$ 
if and only if there is a  $(p,q)$-chain in $R$, or equivalently, the 
inclusion-minimal relation that is transitive and that contains $R$ (see e.g. \cite[p.39]{matouvsek2009invitation}).
It is straightforward to verify that $\tc$ is indeed a closure operator on $S = A\times A$ for all
relations $R$ on $A$.

\paragraph{Graphs.}
A \emph{directed graph $G=(V,E)$} is a tuple with non-empty vertex set $V(G)\coloneqq V$ and arc set $E(G)\subseteq V\times V$.
A \emph{subgraph} of $G$ is a directed graph $H=(V',E')$ such that $V'\subseteq V$ and $E'\subseteq E$.
  We put $\outdeg_G(v)\coloneqq\left|\left\{u\in V \colon (v,u)\in
		E\right\}\right|$ and $\indeg_G(v)\coloneqq\left|\left\{u\in V \colon (u,v)\in
		E\right\}\right|$ to denote the \emph{out-degree} and \emph{in-degree} of a vertex $v$, 
        respectively. A vertex $v$ with $\outdeg_G(v)=0$ is a \emph{leaf} of $G$ and 
		a vertex $v$ with $\indeg_G(v)=0$ is a \emph{root} of $G$. 
        A \emph{hybrid} is a vertex
        $v$ with $\indeg_G(v)\geq 2$. Note that leaves might also be hybrids. 
	   A directed graph $G$ is \emph{phylogenetic} if it does not contain a vertex $v$ such that $\outdeg_G(v)=1$ and
		$\indeg_G(v)\leq 1$, i.e., a vertex with only one outgoing arc and at most one incoming arc. 
	   We sometimes use $u\to v$ to denote the arc $(u,v)\in E(G)$ and 
		$u\leadsto v$ to denote a directed $uv$-path in $G$. 
        An arc $e=(u,w)$ in a DAG $G$ is a
\emph{shortcut} if there is a directed $uw$-path that does not contain the arc $e$
\cite{linz2020caterpillars, DOCKER2019129}. If $G$ has no shortcuts, it is \emph{shortcut-free}.

      Directed graphs $G$ without directed 
      cycles are called \emph{directed acyclic graphs (DAGs)} \cite{book:digraph}. In particular, DAGs do not contain arcs of the form $(v,v)$. Let $G$ be a DAG. 
		We write $v\preceq_G w$ if and only if there is a directed $wv$-path in $G$.  \REV{It is widely known (and indeed easy to verify) that $(V(G),\preceq_G)$ is a poset for every DAG $G$.}
		If $v\preceq_G w$ and $v\neq w$, we write $v\prec_G w$. 
		If $u\to v$ is an arc in $G$, then $v\prec_G u$ and we call $v$ a \emph{child} of $u$ and $u$ a \emph{parent} of $v$. In a DAG $G$ where $u\preceq_G v$, we call $u$ a \emph{descendant} of $v$ and $v$ an \emph{ancestor} of $u$. 
        If $u\preceq_G v$ or $v\preceq_G u$, then $u$ and $v$ are \emph{$\preceq_G$-comparable} and, 
            otherwise, \emph{$\preceq_G$-incomparable}.
            
    For directed graphs $G=(V_G,E_G)$ and $H=(V_H, E_H)$, an
\emph{isomorphism between $G$ and $H$} is a bijective map $\varphi\colon V_G\to V_H$ such that 
$(u,v)\in E_G$ if and only if $(\varphi(u),\varphi(v))\in E_H$. 
If such a map exist, then $G$ and $H$ are \emph{isomorphic}, in symbols $G\simeq H$.

If $G$ is a DAG whose set of leaves is $X$,
then $G$ is a \emph{DAG on $X$}.
Moreover, if a DAG $G$ has a single root, then it is called a {\em network}. 
A \emph{(rooted) tree} is a network that does not contain vertices
with $\indeg_G(v)>1$. A \emph{star-tree} on $X$ is a tree with a single root $\rho$ and
arcs $\rho\to x$ for all $x\in X$.

	\REV{We follow here a notion of phylogenetic networks that has been used, for example, 
    in~\cite{Cardona:09,Hellmuth2023}. However,
    we emphasize that phylogenetic networks are often equipped with the
    additional requirement that every leaf has in-degree one. We do not impose this condition, that is, leaves may have in-degree greater than one and could, therefore, be hybrids as well.
    Since the requirement that leaves have in-degree one is not needed for any of our
	arguments, we adopt this more general definition throughout.
	}

\REV{
For a poset $({Q},\leq)$, the \emph{Hasse diagram}
$\Hasse(Q,\leq)$ is the DAG with vertex set
$Q$ and arcs $(A,B)$ if (i) $B\leq A$ and $A\neq B$ and (ii) there is
no $C\in Q$ with $B\leq C\leq A$ and $C\neq A,B$. 
Since the partial order $\leq$ is transitive and anti-symmetric, $\Hasse(Q,\leq)$ is always a DAG. Moreover, by definition, $\Hasse(Q,\leq)$ is shortcut-free.}

\paragraph{Least common ancestors.}
For a given DAG $G$ on $X$ and a non-empty subset $A\subseteq X$, a vertex $v\in V(G)$ is a \emph{common
ancestor of $A$} if $v$ is an ancestor of every vertex in $A$. Moreover, $v$ is a \emph{least common
ancestor} (LCA) of $A$ if $v$ is a $\preceq_G$-minimal vertex that is an ancestor of all vertices in
$A$. The set $\LCA_G(A)$ comprises all LCAs of $A$ in $G$.

\begin{figure}
  \centering
  \includegraphics[width=0.8\textwidth]{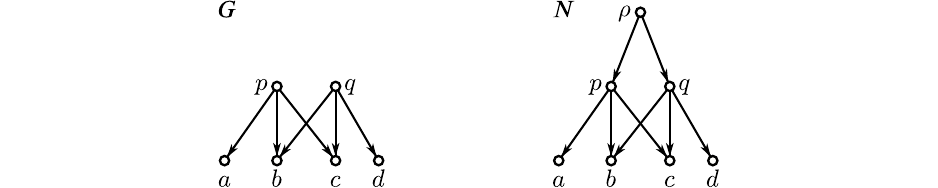}
  \caption{A DAG $G$ and a network $N$, both having leaf-set $X=\{a,b,c,d\}$.
  Note that removing the vertex $\rho$ and its incident arcs from $N$ yields the DAG $G$. 
            In $N$, we have $\LCA_N(\{x,y\}) \neq \emptyset$ for all $x,y\in X$. 
            Moreover, $\LCA_G(\{a,d\}) = \emptyset$ while $\LCA_N(\{a,d\}) = \{\rho\}$ and so $\lca_N(ad)\coloneqq \lca_N(\{a,d\})$ is well-defined whereas $\lca_G(\{a,d\})$ is not.
            In addition, in both DAGs $G$ and $N$ 
            we have $\LCA_G(\{b,c\}) = \LCA_N(\{b,c\})=\{p,q\}$, i.e., the LCA of $b$ and $c$ is not well-defined in either $N$ or $G$. Both $G$ and $N$ are 2-lca-relevant, since e.g.\ $p=\lca_G(ab)=\lca_N(ab)$, $q=\lca_G(cd)=\lca_N(cd)$ and $\rho=\lca_N(ad)$.
            }
  \label{fig:exmpl-LCA}
\end{figure}

In general, not every set $A\subseteq
X$ has a (unique) least common ancestor in a DAG on $X$; for example, consider the DAG $G$
in  Figure~\ref{fig:exmpl-LCA} where $|\LCA_G(\{b,c\})|>1$ and $\LCA_G(\{a,d\})=\emptyset$. 
In a network $N$, the unique root is a common ancestor for all $A\subseteq L(N)$ and, therefore,
$\LCA_N(A)\neq\emptyset$. Moreover,  we are
interested in DAGs where $|\LCA_G(A)|=1$ holds for certain subsets
$A\subseteq X$. For simplicity, we will write
$\lca_G(A)=v$ in case that $\LCA_G(A)=\{v\}$  and say that
\emph{$\lca_G(A)$ is well-defined}; otherwise, we leave $\lca_G(A)$
\emph{undefined}. To recall, we often write $xy$ for sets $\{x,y\}$, which allows
us to put $\lca_G(xy) \coloneqq \lca_G(\{x,y\})$, 
if $\lca_G(\{x,y\})$ is well-defined.  A DAG $G$ on $X$ is {\em 2-lca-relevant} if, 
for all $v\in V(G)$ there are $x,y\in X$ such that $v=\lca_G(xy)$. 
Put simply, in a 2-lca-relevant DAG each vertex is the unique
LCA for one or two leaves. In addition, a DAG $G$ is \emph{lca-relevant}, if
for all $v\in V(G)$ there is a subset $A\subseteq X$ such that $v=\lca_G(A)$. 
Hence, 2-lca-relevant DAGs are lca-relevant.
We note in passing that 
the property of being 2-lca-relevant can be tested in polynomial time; see \cite[Obs.~7.6]{HL:24}.
We illustrate these definitions in Figure~\ref{fig:exmpl-LCA}.

We conclude this section with an important definition that 
relates LCAs in DAGs. As discussed in the 
introduction one of our main aims is to find a DAG that realizes a given set
of LCA constraints. With this in mind we define two key 
relations on LCAs which naturally arise from any DAG $G$.

\begin{definition}[The relations $\srel_G$ and $\rel_G$]\label{def:real_from_DAG}
	For a given DAG $G$ on $X$, we define the relations $\rel_G$ and $\srel_G$ 
    on $\pairs(X)$ as 
    \[\srel_G \coloneqq \{(ab,xy) \mid \lca_G(ab), \lca_G(xy) \text{ are well-defined and } \lca_G(ab)\prec_G \lca_G(xy)\}\] and
	\[\rel_G \coloneqq \{(ab,xy) \mid \lca_G(ab), \lca_G(xy) \text{ are well-defined and } \lca_G(ab)\preceq_G \lca_G(xy)\}.\]
\end{definition}

Note that, by definition, $\srel_G \subseteq \rel_G$.
Moreover, if $G$ is a DAG on $X$, then, $\lca_G(xy) \preceq_G \lca_G(xy)$ for all $x,y \in X$ 
if and only if $\lca_G(xy)$ is well-defined in~$G$, which in turn 
is equivalent to $(xy,xy) \in \rel_G$.  
Hence, $\support_{\rel_G}$ comprises precisely those pairs $xy \in \pairs(X)$ 
for which $\lca_G(xy)$ is well-defined. 
In contrast, $\support_{\srel_G}$ could leave out pairs $xy \in \pairs(X)$ 
for which $\lca_G(xy)$ is well-defined. For example, if $a$ is a leaf in $G$
but not incident to any arc, then $\lca_G(aa)$ is well-defined but 
$\lca_G(aa)\not\prec_G \lca_G(xy)$ for any $x,y\in X$, i.e., 
$aa\notin \support_{\srel_G}$. Nevertheless, as we shall see later 
in Proposition~\ref{prop:cor:G-rel-relplus}, $\support_{\srel_G}^+=\support_{\rel_G}^+$
always holds.

\REV{
To help the reader navigate the different notations used throughout the paper, we provide a collection of the most important definitions in
Table~\ref{tab:sum-defs}. Some of them will be defined more
precisely in the subsequent sections.}

\begin{center}
\begin{table}[ht] \small
  \caption{Summary of main definitions used in this paper.}
  \setlength{\tabcolsep}{8pt} 
  \renewcommand{\arraystretch}{1.4} 
\REV{   \begin{tabular}{p{12.cm}|p{2.5cm}l}	\hline
  For a relation $R$ on $\pairs(X)$ and a DAG $G$ on $X$: & Ref.    \\  \hline
	$\support_R \coloneqq \{p\in \pairs(X) \mid \text{ there is some } q\in \mh{\pairs(X)} \text{ with }  p R q \text{ or } q R p\}$ & Sec.~\ref{sec:basic} \\
	$\support_{R}^+ \coloneqq \support_{R}\ \cup\ \{xx\mid x\in X\} $ & Sec.~\ref{sec:basic} \\ 
	  	 $\srel_G \coloneqq \{(ab,xy) \mid \lca_G(ab), \lca_G(xy) \text{ are well-defined and } \lca_G(ab)\prec_G \lca_G(xy)\}$ & Def.~\ref{def:real_from_DAG} \\
	$\rel_G \coloneqq \{(ab,xy) \mid \lca_G(ab), \lca_G(xy) \text{ are well-defined and } \lca_G(ab)\preceq_G \lca_G(xy) \}$ & Def.~\ref{def:real_from_DAG} \\
	$+$-Closure $R^+$ obtained from $R$ by a finite number of applications of \axiom{R1}, \axiom{R2} and \axiom{R3}  & Def.~\ref{def:rt2c-closure} \& Thm.~\ref{thm:char} \\
	$\cG_R$ is the canonical DAG of $R$ &  Def.~\ref{def:canonG}\\
	$\cN_R$ is the canonical network of $R$ obtained from $\cG_R$ &  Def.~\ref{def:canon-network}\\
	Classical closure $\cl(R) \coloneqq \bigcap_{G\in \mathfrak{G}} \rel_G$ with $\mathfrak{G}$ as the set of DAGs $G$ on $X$ that realize $R$. & Def.~\ref{def:closure}\\
	\hline
	  \end{tabular} }
  \label{tab:sum-defs}
\end{table}
\end{center}


\section{Realizability of LCA Constraints}
\label{sec:realizabilty}

In this section, we define what it means for a DAG to realize a given 
collection of LCA constraints. We consider two variants of realizability: 
one motivated by the approach taken in \cite{Aho:81} and one that naturally generalizes it.

First, recall that, as discussed in the introduction, 
Aho et al.\ \cite{Aho:81} considered the following problem. Given a
collection of constraints of 
the form $(a,b) <(c,d)$, $a \neq b$, $c \neq d$, with $a,b,c,d \in X$, 
can we construct a phylogenetic tree $T$ on $X$ 
so that if $(a,b)<(c,d)$, then $\lca_T(ab) \prec_T \lca_T(cd)$ or determine 
that no such tree exists? This motivates our first
definition for realizability.

\begin{definition}[Strict realization]\label{def:strict-real}
	A   relation $R \subseteq \pairs(X)\times\pairs(X)$ is \emph{strictly realizable} if there is a DAG $G$ on $X$ such that, for all 
	$ab,cd\in\support_R^+$, the vertices $\lca_G(ab)$ and $\lca_G(cd)$ are well-defined and the following implication holds:
	 \begin{itemize}
      \item[\axiom{I0}] $(ab,cd)\in R$ implies that $\lca_G(ab)\prec_G\lca_G(cd)$.
   \end{itemize} 
   \noindent When \axiom{I0} holds, we say that $R$ is \emph{strictly realized} by $G$. 
\end{definition}

We could relax this definition by insisting that
$R$ is ``realizable'' if there is a DAG $G$ such that, for all $ab\ R\ cd$, 
we have $\lca_G(ab) \preceq_G \lca_G(cd)$.
However, with this definition, the star tree (with appropriate leaf set) would \emph{always} realize $R$, 
provided that there are no pairs $ab\ R\ xx$ with $ab \neq xx$.
Hence, this type of ``realizability''  would be too permissive and becomes essentially trivial to satisfy. In order to generalize ``strict realization'', we provide the following definition.

\begin{definition}[Realization]\label{def:real}
	A   relation $R \subseteq \pairs(X)\times\pairs(X)$ is \emph{realizable} if there is a DAG $G$ on $X$ such that, for all 
	$ab,cd\in\support_R^+$ the vertices $\lca_G(ab)$ and $\lca_G(cd)$ are well-defined and the following two implications hold:
       \begin{itemize}
            \item[\axiom{I1}] 
            \REVII{$(ab,cd)\in R$}
            and $(cd,ab)\notin\tc(R)$  implies that $\lca_G(ab)\prec_G\lca_G(cd)$
            \item[\axiom{I2}] 
            \REVII{$(ab,cd)\in R$}
            and $(cd,ab)\in \tc(R)$ implies that  
                         $\lca_G(ab)=\lca_G(cd)$
       \end{itemize}                  
     \noindent When \axiom{I1} and \axiom{I2} hold, we say that $R$ is \emph{realized} by $G$.
\end{definition}

\begin{figure}
  \centering
  \includegraphics[width=0.8\textwidth]{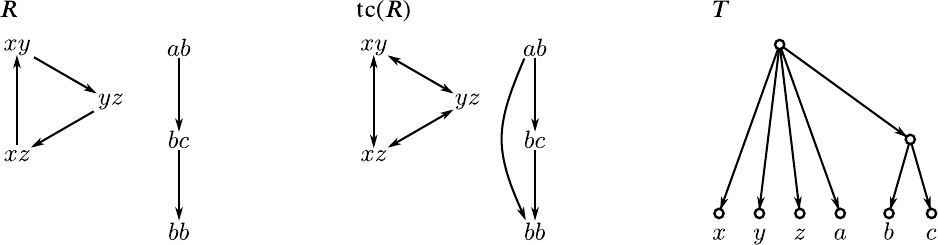}
  \caption{On the left we give a graphical representation of a relation relation $R$ 
  whose vertex set is $\support_R =\{ab,bb,bc,xy,xz,yz\}$. Here we 
  draw an arc $p\to q$ precisely if $q\ R\ p$.
  Similarly, in the middle we give the graphical representation of the transitive closure $\tc(R)$. 
  The phylogenetic tree $T$ on the right realizes $R$.
  }  
  \label{fig:real-1}
\end{figure}

As an example for Definition~\ref{def:real} consider the relation $R$ pictured 
in Figure~\ref{fig:real-1}. In this example, 
\REVII{$(xy,yz)\in R$ and} $(yz,xy)\in \tc(R)$ which, according to \axiom{I2}, implies that 
$\lca_G(xy)=\lca_G(yz)$ must hold for any DAG $G$ realizing $R$. Moreover, 
$(bc,ab)\in \REVII{R}$ and $(ab,bc)\notin \tc(R)$ together with \axiom{I1} implies that 
$\lca_G(bc)\prec_G\lca_G(ab)$ must hold for any DAG $G$ realizing $R$. In this way, 
it is straightforward to verify that the phylogenetic tree $T$ in Figure~\ref{fig:real-1} realizes $R$.

Although the notion of realizability as in Definition~\ref{def:real} looks,  
at a first glance, more restrictive than Definition~\ref{def:strict-real} 
it is in fact a generalization of strict realization, as shown next.

\begin{lemma}\label{lem:str-real=>real}
     Suppose that $R$ is a relation on $\pairs(X)$. 
   Then $R$ is strictly realized by $G$ if and only if $R$ is realized by $G$ and $\tc(R)$ is asymmetric. 
\end{lemma}
\begin{proof}
   \REVII{Let $R$ be a relation on $\pairs(X)$. First assume there is a DAG $G$ on $X$ that strictly realizes $R$. We first show that $\tc(R)$ is asymmetric. For contradiction, assume it is not and let $p,q\in\pairs(X)$ be elements for which $(p,q)\in\tc(R)$ and $(q,p)\in\tc(R)$. By definition of $\tc(R)$ there exists a $(p,q)$-chain in $R$ that can be written as $p = p_0\ R\ p_1\ R\ \cdots\ R\ p_{k-1} R\ p_k=q$, where each $p_i \in \pairs(X)$. Since $G$ strictly realizes $R$, we have 
   $\lca_G(p) = \lca_G(p_0)\prec_G\lca_G(p_1)\prec_G \cdots \prec_G\lca(p_k)=\lca_G(q)$ and, therefore, $\lca_G(p)\prec_G \lca_G(q)$. By analogous arguments, $(q,p)\in\tc(R)$ implies that $\lca_G(q)\prec_G \lca_G(p)$. In summary, $\lca_G(p)\prec_G \lca_G(q)$ and $\lca_G(q)\prec_G \lca_G(p)$ must hold; a contradiction.
Therefore, $\tc(R)$ is asymmetric. We now show that $G$ realizes $R$. By definition, $\lca_G(ab)$ and $\lca_G(cd)$ are well-defined for all $ab,cd\in\support_R^+$.
Let  $(ab,cd)\in R$. Since $G$ strictly realizes $R$, 
 $\lca_G(ab)\prec_G\lca_G(cd)$ must hold. Moreover, 
 $R\subseteq \tc(R)$ and asymmetry of $\tc(R)$ implies that $(cd,ab)\notin\tc(R)$. Hence, $(ab,cd)\in R$ implies that always
 the pre-conditions of \axiom{I1} hold (that is, \axiom{I2} is vacuously satisfied). Since $\lca_G(ab)\prec_G\lca_G(cd)$, \axiom{I1} is indeed always satisfied. 
 Hence, $G$ realizes $R$.}

   \REVII{Suppose now that $R$ is realized by a DAG $G$ and that $\tc(R)$ is asymmetric. In particular, $\lca_G(ab)$ and $\lca_G(cd)$ are well-defined for all $ab,cd\in\support_R^+$. Since for any $ab\ R\ cd$ we also have $(ab,cd)\in\tc(R)$, asymmetry of $\tc(R)$ and \axiom{I1} implies that $\lca_G(ab)\prec_G\lca_G(cd)$. Hence, $G$ indeed satisfies \axiom{I0}, i.e., $G$ strictly realizes $R$.}
\end{proof}

In the rest of this section we shall make some further remarks concerning
our two notions of realizability as given in Definition~\ref{def:strict-real} respectively \ref{def:real} and their relationship. 

In our next result we show that, as expected, for any DAG $G$ 
the relations $\rel_G$ and $\srel_G$ are both realized by $G$.  

\begin{lemma}\label{lem:G-realizes-relG}
    Every DAG $G$ realizes $\rel_G$ and strictly realizes $\srel_G$. 
\end{lemma}
\begin{proof}
    Let $G$ be a DAG on $X$ and define $R\coloneqq \rel_G = \{(ab,cd) \mid \lca_G(ab), \lca_G(cd) \text{
    are well-defined and } \lca_G(ab)\preceq_G \lca_G(cd)\}$. It is easy to verify that $\preceq_G$
    is transitive and that, therefore, $R$ is transitive. Consequently, $\tc(R)=R$. We show now that
    \axiom{I1} and \axiom{I2} are satisfied. For \axiom{I1}, let $(ab,cd)\in R=\tc(R)$ and
    $(cd,ab)\notin\tc(R)=R$. By definition of $R$, we have $\lca_G(ab)\preceq_G \lca_G(cd)$ and
    $\lca_G(cd)\not\preceq_G \lca_G(ab)$. Hence, $\lca_G(ab)\prec_G\lca_G(cd)$ must hold, i.e.,
    \axiom{I1} is satisfied. For \axiom{I2}, suppose $(ab,cd),(cd,ab)\in\tc(R)=R$. By definition of
    $R$, we have $\lca_G(ab)\preceq_G\lca_G(cd)$ and $\lca_G(cd)\preceq_G\lca_G(ab)$. That is,
    $\lca_G(ab)=\lca_G(cd)$ and \axiom{I2} holds.

    To see that $G$ strictly realizes $\srel_G$ we must show that \axiom{I0} holds. But, by
    definition of $\srel_G$, $ab\srel_G cd$ if and only if $\lca_G(ab)$ and $\lca_G(cd)$ are
    well-defined and satisfy $\lca_G(ab)\prec_G\lca_G(cd)$. Thus, \axiom{I0} is trivially satisfied
    and, therefore, $G$ strictly realizes $\srel_G$.    
\end{proof}

Many of the examples considered in this paper are, for simplicity, realizable by trees. 
However, we emphasize that this is not always the case; there exist relations that can be realized 
by some network but not by any tree, see the example in Figure~\ref{fig:intro}.

We now show that if $R$ is realized by a DAG $G$ on $X$, then $R\subseteq\rel_G$.

	\begin{lemma}\label{lem:rel-subset} 
    If $R$ is a relation $\pairs(X)$ that is realized by the DAG $G$ on $X$, 
    then $\lca_G(ab)$ is well-defined for all $ab\in \support^+_R$ and $\lca_G(ab)\preceq_G\lca_G(cd)$ must hold for all $ab\ R\ cd$. 
    Hence, $R\subseteq\rel_G$. 
	\end{lemma} 
    \begin{proof}
    Note that, for all relations $R$, $\support^+_R\setminus \support_R \subseteq \{xx\mid x\in
    X\}$. Hence, for all distinct $a,b\in X$ with $ab\in \support_{R}^+$ we must have $ab\in
    \support_R$. Thus, if $G$ is a DAG on $X$ that realizes $R$, then $\lca_G(ab)$ must be
    well-defined for such $ab\in \support_{R}^+$ and, if $aa\in \support_{R}^+$, then $\lca_G(aa)=a$
    trivially holds. In particular, it is an easy task to verify that, if $G$ realizes $R$, then
    $\lca_G(ab)\preceq_G\lca_G(cd)$ must hold for all $ab\ R\ cd$, irrespective of which implication
    in \axiom{I1} or \axiom{I2} of Definition~\ref{def:real} applies. Hence, $R\subseteq\rel_G$.
    \end{proof}

Note that the converse of Lemma~\ref{lem:rel-subset} is not necessarily true, 
i.e., we may have $R\subseteq \rel_G$ for some DAG $G$ 
even though $G$ does not realize $R$. 
Hence, the notion of realizability given in Definition~\ref{def:real} is 
stronger than just assuming $R\subseteq\rel_G$.
To see this,  consider the following example.
 
\begin{figure}
  \centering
  \includegraphics[width=0.8\textwidth]{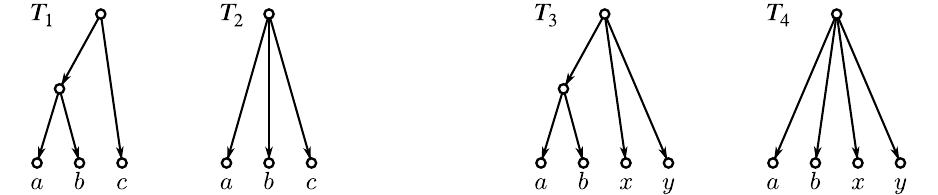}
  \caption{Four phylogenetic trees $T_1$, $T_2$, $T_3$ and $T_4$ used in Example~\ref{empl:subset-non-realizing}, \ref{exmpl:not-subset-real}, \ref{exmpl:R+intro} and \ref{exmpl:union-not-realized}.}
  \label{fig:realize-subset-ex}
\end{figure}

\begin{example}[Subsets of $\rel_G$ that are not realized by $G$]\label{empl:subset-non-realizing}
 Consider the star-tree $T_2$ on $X = \{a,b,c\}$ in Figure~\ref{fig:realize-subset-ex}. We have
 $\support_{\rel_{T_2}} = \{aa,bb,cc,ab,ac,bc\} = \pairs(X)$, as for all pairs in $\pairs(X)$, the respective LCA exists and is unique in $T_2$. 
 Consider now the relation $R = \{(ab,bc)\} \subseteq \rel_{T_2} $. 
 Clearly,
 \REVII{$(ab,cd)\in R$} and $(bc,ab)\notin \tc(R)$. By \axiom{I1}, 
 it must hold that $\lca_G(ab)\prec_G\lca_G(bc)$ for any DAG $G$ realizing $R$.
 Hence, $T_2$ does not realize $R$ since $\lca_{T_2}(ab)=\lca_{T_2}(bc)$.
 Note that $R$ is realized by the tree $T_1$ in Figure~\ref{fig:realize-subset-ex}.
\end{example}

Interestingly, we can even have subsets $R\subseteq \rel_G$ for some DAG $G$
which are not realizable by \emph{any} DAG as the following example shows.

\begin{example}[Subsets of realizable relations that are not realizable]\label{exmpl:not-subset-real}
    Consider the star-tree $T_4$ on $X = \{x,y,a,b\}$, illustrated in Figure~\ref{fig:realize-subset-ex} which, by Lemma~\ref{lem:G-realizes-relG}, realizes $\rel_{T_4}$. 
    Clearly, $\rel_{T_4}$ contains the subset $R =
    \{(xx,ab), (yy,ab), (ab,xy)\}$.  Hence, the relation $R$
    requires, by \axiom{I1}, $\lca_G(ab)\prec_G \lca_G(xy)$ for any DAG $G$ realizing $R$. Since $\lca_{T_4}(ab)= \lca_{T_4}(xy)$, the tree $T_4$ does not realize $R$. Even
    more, $R$ cannot be realized by any DAG $G$. To see this, assume for contradiction, that $G$ is
    a DAG that realizes $R$. By the latter arguments, $\lca_G(ab)\prec_G \lca_G(xy)$ holds. At the
    same time $(xx,ab), (yy,ab)\in R$ together with Lemma~\ref{lem:rel-subset} implies that
    $x=\lca_G(xx)\preceq_G \lca_G(ab)$ and $y=\lca_G(yy)\preceq_G \lca_G(ab)$ holds. As $\lca_G(xy)$
    exists in $G$ and since $\lca_G(ab)$ is a common ancestor of $x$ and $y$ it follows that
    $\lca_G(xy)\preceq_G \lca_G(ab)$ must hold (see also Lemma~\ref{lem:simpleLCAinfer}). This,
    however, violates the requirement that $\lca_G(ab)\prec_G \lca_G(xy)$. Hence, $R\subseteq
    \rel_{T_4}$ is not realizable by any DAG $G$.
\end{example}

As Examples~\ref{empl:subset-non-realizing} and~\ref{exmpl:not-subset-real} show, 
the property that $R \subseteq \rel_G$ for some DAG~$G$ neither implies that $G$ realizes~$R$ 
nor that $R$ is realizable at all.  
In contrast, we have the following result for the relation $ \srel_G$.

\begin{lemma}\label{lem:Rsubset-srelG}
   Let $R$ be a relation on $\pairs(X)$ and $G$ a DAG on $X$. Then, $R$ is strictly realized by $G$ if and only if $R\subseteq\srel_G$.
\end{lemma}
\begin{proof}
    Let $G$ be a DAG on~$X$ and $R$ a relation on $\pairs(X)$. If $R$ is strictly realized by $G$ then, by definition, Condition~\axiom{I0} is satisfied i.e. $\lca_G(ab)$ and $\lca_G(cd)$ are well-defined and satisfy $\lca_G(ab)\prec_G\lca_G(cd)$ for all $(ab,cd)\in R$. By definition of $\srel_G$, we thus have $R\subseteq\srel_G$.
    Conversely, suppose $R \subseteq \srel_G$. 
    By definition of $\srel_G$ and since $R \subseteq \srel_G$, for all $(ab,cd) \in R$, we have $\lca_G(ab) \prec_G \lca_G(cd)$.  
    Hence, \axiom{I0} is satisfied. 
    By definition, $R$ is therefore strictly realized by~$G$. 
\end{proof}


\section{The \(+\)-Closure}
\label{sec:closure}

Our next goal is to characterize when a relation $R$ on $\pairs(X)$  
can be realized by some DAG (or network) on $X$. 
To this end, we define and study in this section a closure relation $R^+$ 
associated with $R$, which contains additional information about any DAG $G$ 
that could potentially realize $R$. 
A related approach for realizing LCA constraints using phylogenetic trees 
is outlined in~\cite[Section~2]{ng1996reconstruction}. 

The key observation underlying the definition of $R^+$ is based on the following result, 
which shows how ancestor relations among LCAs can be inferred from other ones.

\begin{lemma}\label{lem:simpleLCAinfer}
Let $G$ be a DAG on $X$ and suppose that $a,b,c,d,x,y \in X$ (not necessarily 
distinct) are such that $\lca_G(ab), \lca_G(ac), \lca_G(xy)$ and $\lca_G(bd)$ 
are all well-defined. Then, the following statement holds: 
\begin{equation}\label{eq:LCA-implication}
\text{if } \lca_G(ac)\ \preceq_G\ \lca_G(xy) \text{ and } \lca_G(bd)\ \preceq_G\ \lca_G(xy), \text{ then }
 \lca_G(ab)\ \preceq_G\ \lca_G(xy). 
 \end{equation}
\end{lemma}
\begin{proof}
Let $G$ be as defined in the statement of the lemma. 
Suppose that  $\lca_G(ac)\ \preceq_G\ \lca_G(xy)$ and
$\lca_G(bd)\ \preceq_G\ \lca_G(xy)$. 
Hence,  $\lca_G(xy)$ is a common ancestor of $a$ and $b$. 
This and the defining property of the least common ancestor $\lca_G(ab)$ implies that
$\lca_G(ab)\ \preceq_G\ \lca_G(xy)$.
\end{proof}

Although the implication in Equation~\eqref{eq:LCA-implication} from 
Lemma~\ref{lem:simpleLCAinfer} is only one of many that can be derived from ``LCA-$\prec_G$ relationships'',  it turns out that this particular implication is sufficient to characterize 
realizable relations and motivates the following definition.

\begin{definition}[Cross-Consistency]\label{def:2cons}
A relation $R$ on $\pairs(X)$ is \emph{cross-consistent} if 
for all $a,b,c,d,x,y \in X$ (not necessarily distinct) 
the following statement holds: 
\[
\mbox{if $ac\ R\ xy$, $bd\ R\ xy$ and $ab\in\support_R$, then $ab\ R\ xy$}. 
\]
\end{definition}

A given relation $R$ on $\pairs(X)$ typically contains only ``partial information'' on pairs of LCAs in $\support_{R}$. 
We now define the aforementioned closure $R^+$ of $R$ that allows us to extend $R$ to encompass information about additional pairs.

\begin{definition}[$+$-Closure]\label{def:rt2c-closure}
 If $R$ is a relation on $\pairs(X)$, then the \emph{reflexive transitive cross-consistent closure (for short, $+$-closure)} is defined as
\[
R^+ \coloneqq \bigcap_{R'\in\mathfrak{R}} R',
\]
where $\mathfrak{R}$ denotes the set of all relations $R'$ on $\pairs(X)$
that are $\support_R^+$-reflexive, transitive, cross-consistent, and satisfy $R \subseteq R'$. 
\end{definition}

It is an easy task to verify that the $+$-closure $R^+$ of a relation $R$ is $\support_R^+$-reflexive, transitive and cross-consistent since $R^+\subseteq R'$ and, therefore, 
 $\support_R^+ \subseteq  \support_{R^+}\subseteq \support_{R'}$
for all $R'\in \mathfrak{R}$. We now show that $R^+$ also satisfies that classical closure axioms. 

\begin{proposition}\label{lem:closure-axioms}
Let $R$ be a relation on $\pairs(X)$.
The operator $^+ \colon R\, \mapsto\, R^+$ is a closure operator on $\pairs(X)\times \pairs(X)$, i.e., it satisfies
\emph{extensivity} $R \subseteq R^+$; 
\emph{monotonicity}  $R_1 \subseteq R_2$ $\implies$ $R_1^+ \subseteq R_2^+$;
and \emph{idempotency} $(R^+)^+ = R^+$.
\end{proposition}
\begin{proof}
Let $R$ be a relation on $\pairs(X)$. Let $\mathfrak{R}$ be 
the set of all relations on $\pairs(X)$ that are $\support_R^+$-reflexive, transitive,
cross-consistent and contain $R$, so that $R^+=\bigcap_{R'\in\mathfrak{R}}R'$.
By definition, every $R'\in\mathfrak{R}$ satisfies $R\subseteq R'$. Hence
$R \subseteq \bigcap_{R'\in\mathfrak{R}} R' = R^+$, i.e., extensivity holds.

Assume that $R_1\subseteq R_2$ for two relations $R_1$ and $R_2$ on $\pairs(X)$. Let
$\mathfrak{R}_{i}$ be the set of all relations on $\pairs(X)$ that are
$\support_{R_i}^+$-reflexive, transitive, cross-consistent and contain $R_i$ for $i\in\{1,2\}$.
Note that $R_1\subseteq R_2$ implies $\support_{R_1}\subseteq\support_{R_2}$ and since both
relations are on $\pairs(X)$ we thus have $\support_{R_1}^+\subseteq\support_{R_2}^+$. Consequently,
every $\support_{R_2}^+$-reflexive relation is, in particular, also $\support_{R_1}^+$-reflexive.
Now, by definition and the fact that $R_1\subseteq R_2$, it follows that any relation
$S\in\mathfrak{R}_2$ is $\support_{R_1}^+$-reflexive, transitive, cross-consistent and contains
$R_1$, so $S\in\mathfrak{R}_1$. In other words, $\mathfrak{R}_2\subseteq\mathfrak{R}_1$.
Consequently,
\[R_1^+=\bigcap_{S\in\mathfrak{R}_1} S\subseteq\bigcap_{S\in\mathfrak{R}_2} S=R_2^+\]
follows, i.e. monotonicity holds.
 
Now let $\mathfrak{R}_{R^+}$ be the set of $\support_{R^+}^+$-reflexive, transitive and cross-consistent relations
on $\pairs(X)$ that contain $R^+$.
By construction $R^+$
itself belongs to $\mathfrak{R}_{R^+}$, hence
\[
(R^+)^+ \;=\; \bigcap_{R'\in\mathfrak{R}_{R^+}} R' \;\subseteq\; R^+.
\]
Moreover, by extensivity, we have $R^+\subseteq (R^+)^+$.
Thus,  $(R^+)^+=R^+$ which proves idempotency.
\end{proof}

We note in passing that for all $ab\in \support_{R}^+$ (where $a=b$ might be possible) 
we have $ab\ R^+\ ab$ since $R^+$ is $\support_R^+$-reflexive. This together with 
cross-consistency of $R^+$ implies \REV{the following:
\begin{observation}\label{obs:always_xx<xy}
    If $ab\in \support_R^+$, then $(aa,ab),(bb,ab) \in R^+$ always holds.
\end{observation}
}

We show now that the relation $R^+$ can be obtained by application of three simple rules on $R$
and can therefore be computed in polynomial time.

\begin{theorem}\label{thm:rules}
Let $R$ be a relation on $\pairs(X)$. 
Assume that $\altsrel$ is a   relation obtained from $R$ by starting with $\altsrel=R$ and 
first applying the rule 
\begin{enumerate}
  \item[\axiom{R1}] \emph{Reflexivity:} for all  $p\in\support^+_R$, add $p\ \altsrel\ p$.
\end{enumerate}
Afterwards,  repeatedly apply one of the following two rules in any order, until none of the rules can be applied:
\begin{enumerate}
  \item[\axiom{R2}] \emph{Transitivity:} if $p\ \altsrel\ q$ and $q\ \altsrel\ r$, add $p\ \altsrel\ r$.
  \item[\axiom{R3}] \emph{Cross-Consistency:} if $ab\in \support_S$ and $ac\ \altsrel\ xy$ and $bd\ \altsrel\ xy$ for some $c,d\in X$,
        then add $ab\ \altsrel\ xy$.
\end{enumerate}
Then, $\altsrel=R^+$ and $\support_R^+ =  \support_{R^+} = \support_S$.
\end{theorem}
\begin{proof}
  Suppose that $\altsrel$ is obtained as outlined in the statement and let 
  $\mathfrak{R}$ be the set of all relations $R'$ on $\pairs(X)$ that contain $R$ and are $\support_R^+$-reflexive, transitive and cross-consistent. 
  Since $R^+ = \cap_{R'\in\mathfrak{R}} R'$ have 
  \begin{equation}\label{eq:subset}
 		R^+ \subseteq R' \text{ and, thus, }   \support_{R^+}\subseteq \support_{R'}
 		 \text{ for all } R'\in\mathfrak{R}.  
  \end{equation}
  Since the rules \axiom{R1}, \axiom{R2} and \axiom{R3} are exhaustively applied, the final relation
  $\altsrel$ must be $\support_R^+$-reflexive, transitive, cross-consistent and satisfies, by definition,
  $R\subseteq \altsrel$. Therefore, $\altsrel\in\mathfrak{R}$ holds, and consequently $R^+\subseteq S$. Moreover, 
  $\support_S=\support_R^+$ holds by construction.

  We now show that $\altsrel\subseteq R^+$ must hold by considering the intermediate
  applications of the rules \axiom{R1}--\axiom{R3}. Specifically, let $\altsrel_0, \altsrel_1,\ldots,\altsrel_n$
  denote the relations on $\pairs(X)$ for which $\altsrel_0=R$, $\altsrel_n=\altsrel$ and where,
  for each $0\leq i\leq n-1$, $\altsrel_{i+1}$ is constructed from $\altsrel_{i}$ by applying
  exactly one of \axiom{R1}, \axiom{R2} and \axiom{R3}. Observe that Proposition~\ref{lem:closure-axioms} ensures that
  $\altsrel_0=R$ must be a subset of $R^+$. Moreover, $\altsrel_1$ is obtained from
  $\altsrel_0=R$ by applying \axiom{R1} and thus $\altsrel_1=R\cup\{(p,p)\mid
  p\in\support_R^+\}$. The latter equality taken together with $R\subseteq R^+$ and the fact
  that $R^+$ is $\support_R^+$-reflexive, imply that $\altsrel_1\subseteq R^+$
  and, thus, $\support_{\altsrel_1}\subseteq \support_{R^+}$.
  
  Now,
  assume that $\altsrel_i$ is contained in $R^+$ for some fixed $i$ with $1\leq i<n$. Thus
  $\altsrel_{i+1}$ is constructed from $\altsrel_{i}$ by applying exactly one of \axiom{R2} and \axiom{R3}. We now consider what happens in each of these two cases.
  
  If $\altsrel_{i+1}$ is obtained from $\altsrel_i$ by applying \axiom{R2}, then
  $\altsrel_{i+1}=\altsrel_i\cup\{(p,r)\}$ for some pair $(p,r)$ for which
  $(p,q),(q,r)\in\altsrel_{i}\subseteq R^+$. By Eq.~\eqref{eq:subset}, 
  $(p,q),(q,r)\in R'$ for all $R'\in \mathfrak{R}$. Moreover,
  since each relation in $\mathfrak{R}$ is transitive, the pair $(p,r)$ 
  must be contained in $R'$ for all $R'\in \mathfrak{R}$. 
  By definition of $R^+$, we have $p\ R^+\ r$.
  Consequently, $\altsrel_{i+1}\subseteq R^+$.

  If $\altsrel_{i+1}$ is obtained from $\altsrel_i$ by applying \axiom{R3}, then
  $\altsrel_{i+1}=\altsrel_i\cup\{(ab,xy)\}$ for some pair $(ab,xy)$ for which
  $ab\in \support_{S_i}$. 
  In particular, \axiom{R3} requires the existence of some $(ac,xy),(bd,xy)\in\altsrel_{i}$.
   Since, by assumption,   $\altsrel_i\subseteq R^+$,   we have
   $(ac,xy),(bd,xy)\in R^+$. By Eq.~\eqref{eq:subset}, 
  $(ac,xy),(bd,xy)\in R'$ and $ab\in\support_{R'}$ for all $R'\in \mathfrak{R}$. Moreover, since each 
  relation in $\mathfrak{R}$ is cross-consistent, the pair $(ab,xy)$ must also be contained
  in $R'$ for all $R'\in \mathfrak{R}$. By definition of $R^+$, we have
  $ab\ R^+\ xy$. Consequently, $\altsrel_{i+1}\subseteq R^+$.

  In both of these cases, $\altsrel_{i+1}\subseteq R^+$ holds and so, 
  by induction, it follows that  $\altsrel=\altsrel_n\subseteq R^+$.
  This together with   $R^+\subseteq \altsrel$ implies that   $R^+ = \altsrel$.
 In particular, $\support_R^+ = \support_S =  \support_{R^+}$.
\end{proof}

\REV{The main steps to compute the $+$-closure $R^+$ of a relation $R$ 
 according to Theorem~\ref{thm:rules} are summarized in Algorithm~\ref{alg:Rplus}.}

\begin{algorithm}
  \caption{\REV{\textsc{Computing the $+$-closure of a Relation $R$}}}
  \label{alg:Rplus}
  \REV{
  \begin{algorithmic}[1]
  \Require  A binary relation  $R$ on $\pairs(X)$
    \Ensure  Returns the $+$-closure $R^+$ of $R$
    \State Compute $\support_{R}^+$ 
    \State $S\gets R$ 
    \ForAll{$p\in \support_{R}^+$} \Comment{Rule \axiom{R1}}
        \State $S\gets S\cup \{(p,p)\}$ 
     \EndFor   
     \State rule\_applicable $\gets$ \texttt{true} 
     \While{rule\_applicable}
         \If{there are $(p,q),(q,r)\in S$ for which $(p,r)\notin S$} \Comment{Rule \axiom{R2}}
            \State $S\gets S\cup \{(p,r)\}$ 
          \ElsIf{there are  $(ac,xy),(bd,xy)\in S$ for which $ab\in \support_S$ and $(ab,xy)\notin S$}  \Comment{Rule \axiom{R3}}
            \State $S\gets S\cup \{(ab,xy)\}$ 
           \Else\  rule\_applicable $\gets$ \texttt{false} 
          \EndIf  
      \EndWhile  
	 	 \State \Return $S$
  \end{algorithmic}}
\end{algorithm}

\REV{
\begin{theorem}\label{thm:R+-polytime}
Algorithm~\ref{alg:Rplus} correctly determines the $+$-closure $R^+$ of a 
relation $R$ on $\pairs(X)$ in polynomial time in $|X|$. 
\end{theorem}}
\begin{proof}
	\REV{
    Let $R$ be a relation on $\pairs(X)$. 	
	By Theorem~\ref{thm:rules}}, $R^+$ can be
  obtained by repeatedly applying Rules \axiom{R1}, \axiom{R2}, and \axiom{R3}, starting with $S =
  R$ and extending it step by step. \REV{It is now an easy task to verify that 
  Algorithm~\ref{alg:Rplus} correctly determines $R^+$.}

  \REV{For the runtime, observe that}
  $S \subseteq \pairs(X) \times \pairs(X)$ and $|\pairs(X)| 
  = \binom{|X|}{2} + |X| \in O(|X|^2)$. Hence, $|S| \in O(|X|^4)$ at each step. Thus,
  checking whether one of the Rules \axiom{R1}, \axiom{R2}, or \axiom{R3} can be applied to $S$ can
  be done in polynomial time, and adding a pair $(p,q)$ to $S$ requires only constant time. Clearly,
  the process of constructing $R^+$ by applying Rules \axiom{R1}, \axiom{R2}, and \axiom{R3} must
  terminate, since at most $O(|X|^4)$ pairs can be added to $S$. In summary, $R^+$ can be
  constructed in polynomial time in $|X|$.
\end{proof}

\begin{figure}
    \centering
    \includegraphics[width=0.8\textwidth]{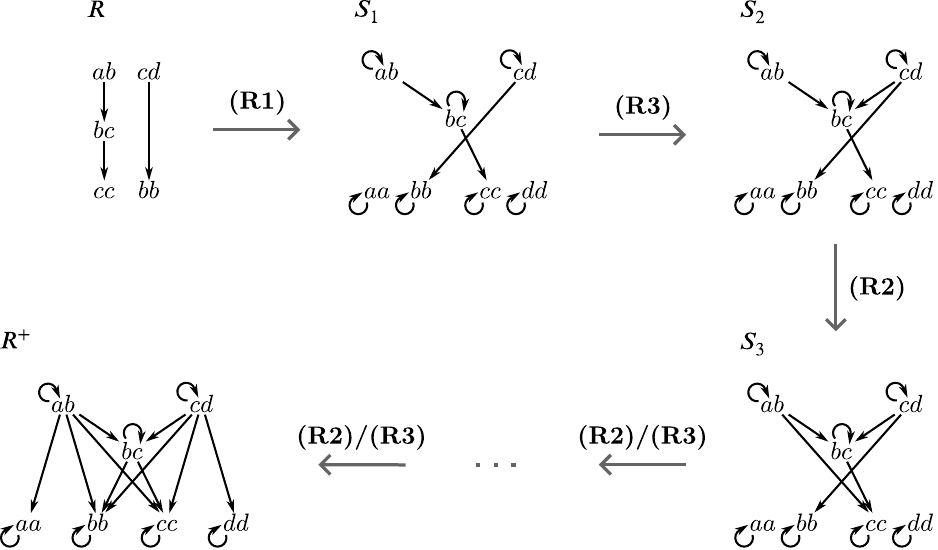}
    \caption{\REV{Shown are several graphical representation of relations $S\in \{R,S_1,S_2,S_3,R^+\}$
            where we draw an arc $p\to q$ precisely if $q\ S\ p$.
            The initial relation is $R= \{(bc,ab), (cc,bc),  (bb,cd)\}$ on $\pairs(X)$ with $X=\{a,b,c,d\}$. 
            Application of \axiom{R1} to $R$ yields the relation $S_1$. 
            Application of \axiom{R3} to $(bb,cd),(cd,cd)\in S_1$ yields the relation $S_2 = S_1\cup \{(bc,cd)\}$. Application of \axiom{R2} to $(cc,bc),(bc,ab)\in S_2$
            yields $S_3 = S_2\cup \{(cc,ab)\}$. Further applications of \axiom{R2}
            and \axiom{R3} result in $R^+$.}
            }
    \label{fig:stepwise-closure}
\end{figure}

\REV{Illustrative examples of relations $R$ and their $+$-closure $R^+$ are
provided in Figure~\ref{fig:stepwise-closure} and \ref{fig:working-example}.}
In the following example, we demonstrate why $R^+$ is useful.
\begin{example}\label{exmpl:R+intro} 
Consider the relation $R$ with $xx\ R\ ab$, $yy\ R\ ab$, $aa\ R\ xy$, and $bb\ R\ xy$,
which is realized by the star-tree $T_4$ on $X = \{a,b,x,y\}$, as illustrated in Figure~\ref{fig:realize-subset-ex}.
Observe first the accordance with Lemma~\ref{lem:simpleLCAinfer}: we have
$x = \lca_{T_4}(xx) \preceq_{T_4} \lca_{T_4}(ab)$ and $y = \lca_{T_4}(yy) \preceq_{T_4} \lca_{T_4}(ab)$, and moreover
$\lca_{T_4}(xy) \preceq_{T_4} \lca_{T_4}(ab)$. In fact, $\lca_{T_4}(xy) = \lca_{T_4}(ab)$ holds.

However, the relation $R$ is not cross-consistent, since $xx\ R\ ab$ and $yy\ R\ ab$ hold, but $xy\ R\ ab$ does not.
In contrast, we have both $xy\ R^+\ ab$ and $ab\ R^+\ xy$. As we shall see later in Lemma~\ref{lem:soundness}, this implies that
$\lca_G(xy) \preceq_G \lca_G(ab)$ and $\lca_G(ab) \preceq_G \lca_G(xy)$,
and thus $\lca_G(xy) = \lca_G(ab)$ must hold for any DAG $G$ realizing $R$. 
\end{example}

Example~\ref{exmpl:R+intro} shows that it is helpful to think of $R^+$ as the \emph{minimal} extension of $R$ for which
$\preceq_G$-comparability requirements between unique LCAs are enforced across all DAGs $G$ that realize $R$.

We conclude this section with some observations and 
results concerning the closure $R^+$ that will be useful later.
Since $R^+$ is transitive, we have $\tc(R^+)= R^+$. This together 
with $R\subseteq R^+$ and extensivity of the transitive closure implies: 

\begin{observation}\label{obs:tcR-subset-R+}
    For all relations $R$ we have $\tc(R)\subseteq R^+$.
\end{observation}

The following result shows that the closure $R^+$ contains valuable 
additional information about any DAG $G$ that can potentially realize $R$.
Indeed, we now show that if a DAG $G$ realizes $R$, then we 
not only have $R\subseteq \rel_G$ but we must have $R^+\subseteq \rel_G$.

\begin{lemma}
\label{lem:soundness}
Assume a relation $R$ on $\pairs(X)$ is realized by a DAG $G$ on $X$.
Then for all $a,b,x,y\in X$
\[
ab\ R^+\ xy \quad\Longrightarrow\quad \lca_G(ab)\ \preceq_G\ \lca_G(xy).
\]
In particular, $R^+\subseteq \rel_G$ for every DAG $G$ that realizes $R$.
\end{lemma}
\begin{proof}
Suppose that $R$ on $\pairs(X)$ is a relation that is realized by the DAG $G$ on $X$. Let
$S_0, S_1,\ldots,S_n$ denote relations on $\support^+_R$ for which $S_0=R$,
$S_n=R^+$ and where, for each $S_{i+1}$ is constructed from $S_{i}$, $0\leq i\leq
n-1$, by applying exactly one of \axiom{R1}, \axiom{R2} and \axiom{R3} to $S_i$; the existence of at least one such a sequence is ensured by Theorem~\ref{thm:rules}.
We will prove that, for each $0\leq i\leq n$, we have
\begin{equation}\label{eq:ind-impl}
ab\ S_i\ xy\implies \lca_G(ab)\ \preceq_G\ \lca_G(xy).
\end{equation}
By Lemma~\ref{lem:rel-subset}, $ab\ R\ xy$ always implies that  
$\lca_G(ab)\ \preceq_G\ \lca_G(xy)$.
Thus, for $i=0$ where $S_i=R$, Eq.~\ref{eq:ind-impl} trivially holds. 
Moreover, by definition, $S_1$ is obtained from $R$ by applying \axiom{R1} and hence $S_1=R\cup\{(p,p)\mid p\in\support_{R}^+\}$. 
Clearly the elements $(p,p)\in S_1\setminus R$ satisfy Eq.~\ref{eq:ind-impl}, since $\lca_G(ab)\preceq_G\lca_G(ab)$ holds for all $p=ab\in\support_{R}^+$.
Therefore, Eq.~\ref{eq:ind-impl} is satisfied for all elements in $S_1$ and thus, for $i\in \{0,1\}$.

Assume next, that Eq.~\ref{eq:ind-impl} holds for some fixed $i\geq 1$ and consider $S_{i+1}$.
By definition, $S_{i+1}=S_i\cup\{(p,q)\}$ for some $p,q\in\support^+_R$. 
Moreover, by the induction hypothesis, Eq.~\ref{eq:ind-impl} holds for all pairs in $S_i = S_{i+1}\setminus \{(p,q)\}$. 
Hence, it suffices to show that Eq.~\ref{eq:ind-impl} hold for $p\ S_{i+1}\ q$. To this end, we
consider now the two possible cases for how $p\ S_{i+1}\ q$ was derived.

If $p\ S_{i+1}\ q$ was added via \axiom{R2}, then there must exist $ab, xy, cd\in\support^+_R$ for which
$p=ab\ S_i\ xy$ and $xy\ S_i\ cd = q$. By Eq.~\ref{eq:ind-impl}, \(\lca_G(ab)\ \preceq_G\
\lca_G(xy)\) and \(\lca_G(xy)\ \preceq_G\ \lca_G(cd)\). Thus, by transitivity of $\preceq_G$
(concatenation of directed paths in $G$), $\lca_G(ab)\preceq_G \lca_G(cd)$ holds as well.

If $p\ S_{i+1}\ q$ was added via \axiom{R3}, then there must exist $ab, ac, bd, xy\in\support^+_R$ for which
$ac\ S_i\ xy$, $bd\ S_i\ xy$ and where $p=ab$, respectively $q=xy$. 
Since $G$ realizes $R$ and $ab, ac, bd, xy\in\support^+_R$  it follows that 
$\lca_G(ab), \lca_G(ac), \lca_G(xy)$ and $\lca_G(bd)$ must be well-defined in $G$. 
By Eq.~\ref{eq:ind-impl},
$\lca_G(ac)\ \preceq_G\ \lca_G(xy)$ and $\lca_G(bd)\ \preceq_G\ \lca_G(xy)$. By
Lemma~\ref{lem:simpleLCAinfer}, $\lca_G(ab)\ \preceq_G\ \lca_G(xy)$ holds. 

In summary, for each of the relations $S_i$, and thus, in particular, for $S_n = R^+$, 
Eq.~\ref{eq:ind-impl} is satisfied. Hence, $ab\ R^+\ xy \quad\Longrightarrow\quad \lca_G(ab)\ \preceq_G\ \lca_G(xy)$
for every DAG $G$ that realizes $R$. This and Def.~\ref{def:real_from_DAG}  implies that 
 $R^+\subseteq \rel_G$ for every DAG $G$ that realizes $R$.
\end{proof}

The final result of this section shows that for every DAG $G$, the underlying relation $\rel_G$
and its $+$-closure agree.

\begin{proposition}\label{prop:G-rel-relplus} 
    For all DAGs $G$ we have $\rel_G= \rel_G^+$.
\end{proposition}
\begin{proof}
    Let $G$ be a DAG on $X$. To show that $\rel_G=\rel_G^+$ it suffices, by Theorem~\ref{thm:rules},
     to show that $\rel_G$ is
    $\support_{\rel_G}^+$-reflexive, transitive and cross-consistent.
    Since $x=\lca_G(xx)$ is well-defined for all $x\in X$
    and $x\preceq_G x$,  we have $xx\rel_G xx$ for all $x\in X$. In particular, it follows that 
    $\support_{\rel_G}=\support_{\rel_G}^+$.  Since  $\preceq_G$ is reflexive, i.e. $\lca_G(ab)\preceq_G\lca_G(ab)$ for all $ab\in\pairs(X)$
     where $\lca_G(ab)$ is well-defined, it follows that $\rel_G$ is reflexive
    and, thus, in particular $\support_{\rel_G}^+$-reflexive. Similarly, transitivity of $\preceq_G$  implies that $\rel_G$ is transitive.
    We conclude the proof of the proposition by showing that $\rel_G$ is cross-consistent. To this end, suppose there are 
    $a,b,c,d,x,y \in X$ (not necessarily distinct) such that $ac \rel_G xy$ and $bd \rel xy$. 
    By definition, $\lca_G(ac) \preceq_G \lca_G(xy)$ and 
    $\lca_G(bd) \preceq_G \lca_G(xy)$ holds. Suppose now that $ab \in \support_{\rel_G}$
    and thus, that  $\lca_G(ab)$ is well-defined. 
    Then we can apply Lemma~\ref{lem:simpleLCAinfer}, 
    to conclude that $\lca_G(ab)\ \preceq_G\ \lca_G(xy)$ holds. 
    Hence, $ab \rel_G xy$ holds and $\rel_G$ is cross-consistent.
    In summary, $\rel_G=\rel_G^+$ holds for all DAGs $G$. 
\end{proof}

As we shall see later in Proposition~\ref{prop:cor:G-rel-relplus}, $\srel_G\neq  \srel_G^+$
for all DAGs $G$.


\section{The Canonical DAG and  Characterization of Realizability}
\label{sec:canonDAG}

In this section, we shall give a characterization of when 
a relation $R$ on $\pairs(X)$ is realizable by a DAG.
Our approach proceeds as follows.
We employ the closure $R^+$ to define an equivalence relation $\sim_{R^+}$ comprising pairs $(p,q)$ that satisfy $p\ R^+\ q$ and $q\ R^+\ p$. Based on $\sim_{R^+}$, 
 we will define the \emph{canonical DAG} $\cG_R$
whose vertices are the equivalence classes of $\sim_{R^+}$  and whose arcs are induced by the partial order on these classes. The classes $[aa]$ correspond to leaves of $\cG_R$ 
and are relabeled by $a$ to obtain a DAG on~$X$. 
We then show in Theorem~\ref{thm:char} that $R$ can be realized by a DAG if and only if it is realized by $\cG_R$, and that this is equivalent to $R$ satisfying the following two simple conditions.

\begin{definition}\label{def:X12}
For a relation $R$ on $\pairs(X)$, we define the following two conditions:
\begin{description}
\item[$\qquad$\textnormal{Condition} \axiom{X1}:]  For all $a,b,x\in X$: $ab\neq xx$ implies \REVII{$(ab, xx)\notin R$.}
\item[$\qquad$\textnormal{Condition} \axiom{X2}:]  For all $a,b,x,y\in X$: \REVII{$(ab, xy)\in R$} and  $(xy, ab)\not\in \tc(R)$ implies $(xy,ab)\notin R^+$. 
\end{description}
\end{definition}

\REVII{
The following result shows, in particular, that if 
$(ab, xx)\notin R$ then even $(ab, xx)\notin R^+$.

\begin{lemma}\label{lemma:X1-R-R+}
    A relation $R$ on $\pairs(X)$ satisfies \axiom{X1} if and only if $R^+$ satisfies \axiom{X1}.  
\end{lemma}
\begin{proof}
Let $R$ be a relation on $\pairs(X)$. Clearly, $R\subseteq R^+$ implies that
 $R$  satisfies \axiom{X1} whenever $R^+$ does. 
 
Now suppose that $R$ satisfies \axiom{X1}. 
By Theorem~\ref{thm:rules}, there exists a sequence $R = S_0, S_1, \dots, S_n = R^+$, $n\geq 0$ 
such that, $S_1$ is obtained from $S_0$ by applying \axiom{R1} and
$S_i$ is obtained from $S_{i-1}$ by application of \axiom{R2} or \axiom{R3} for each 
$i \in \{2, \dots, n\}$. We show now that $S_i$ satisfies \axiom{X1} for all $i \in \{0,1, \dots, n\}$. 
Clearly, $S_0=R$ satisfies \axiom{X1}. Moreover, $S_1\setminus S_0\subseteq\{(xx,xx)\mid x\in X\}$ by definition of \axiom{R1}, so $S_1$ also satisfies \axiom{X1}. Thus we can, by induction, assume that 
$S_k$ satisfies \axiom{X1} some $1\leq k < n$. 
Consider now $S_{k+1}$ obtained from $S_k$ by application of either \axiom{R2} or \axiom{R3}. 
Since $S_k$ satisfies \axiom{X1}, 
it holds that $a,b,x \in X$ with $ab \neq xx$ implies $(ab,xx) \notin S_k$. 

Suppose first that we applied rule \axiom{R2} to obtain $S_{k+1}$ from $S_k$. 
In this case, there exists $(p,q), (q,r) \in S_{k}$ and we construct $S_{k+1} = S_k\cup \{(p,r)\}$.
By assumption, neither $(p,q)$ nor $(q,r)$ is of the form $(ab,xx)$ with $ab \neq xx$. 
Hence it is straightforward to verify that,  in this case, $(p,r)$ cannot be of the 
the form $(ab,xx)$ with $ab \neq xx$. This and the fact that  $S_{k}$  satisfies \axiom{X1} 
implies that $S_{k+1}$ satisfies \axiom{X1}.  

Finally suppose that $S_{k+1}$ is obtained from $S_{k}$ by application of \axiom{R3}. 
Assume, for contradiction, that $(ab,xx)\in S_{k+1}$ for some 
$a,b,x \in X$ with $ab \neq xx$. As $S_k$ satisfies \axiom{X1}, it follows that $(ab,xx) \notin S_k$. 
Since $S_k$ and $S_{k+1}$ differ only by one element, it follows that $S_{k+1} = S_k\cup \{(ab,xx)\}$.
Since \axiom{R3} was applied to obtain $S_{k+1}$, 
there must be elements $(ac, xx), (bd, xx)\in S_k$ such that $ab\in\support_{S_k}$. 
However, since $S_k$ satisfies \axiom{X1} and $(ac, xx),(bd, xx) \in S_k$
it follows that $ac=xx$ and $bd=xx$. Hence,  $ab = xx$; a contradiction. 
Thus, after application of \axiom{R3}, $S_{k+1}$ cannot contain elements
$(ab,xx)$ with  $ab \neq xx$. Therefore, $S_{k+1}$ satisfies \axiom{X1}. 

In summary and by induction, $S_n=R^+$ satisfies \axiom{X1}. 
\end{proof}
}

We continue by showing that \axiom{X1} and \axiom{X2} are necessary conditions 
for a relation $R$ to be realizable. 

\begin{lemma}\label{lem:no-collapse}
If $R$ on $\pairs(X)$ is realizable, 
then $R$ satisfies \axiom{X1} and \axiom{X2}. 
\end{lemma}
\begin{proof}
Suppose that $R$ on $\pairs(X)$ is realized by the DAG $G$ on $X$. 
\REVII{We show the contrapositive of \axiom{X1}. To that end, assume $a,b,x\in X$ are elements such that 
$(ab,xx)\in R$.}
By \REVII{Lemma~\ref{lem:rel-subset}}, $\lca_G(ab)\ \preceq_G\ \lca_G(xx)=x$.
Since $x$ is a leaf of $G$, it does not have any children
and thus, $\lca_G(ab) = \lca_G(xx)=x$ which immediately implies that $a=b=x$. 
\REVII{In other words, $ab=xx$.}

For Condition \axiom{X2}, suppose that 
\REVII{$(ab, xy)\in R$ and}  $(xy, ab)\notin \tc(R)$.
\REVII{Since $G$ realizes $R$, \axiom{I1} ensures that $\lca_G(ab)\prec_G\lca_G(xy)$ and, thus, that $\lca_G(xy)\not\preceq_G\lca_G(ab)$. The latter together with Lemma~\ref{lem:soundness} implies that $(xy,ab)\notin R^+$ and we can conclude that \axiom{X2} holds.} 
\end{proof}

    \REV{The following definition uses the standard construction that associates a poset with
    a reflexive and transitive relation (also known as a \emph{preorder} or \emph{quasi-order}),
    by identifying precisely those elements that violate anti-symmetry.
    In our setting, this construction is not applied directly to $R$, which need not be reflexive or transitive, but rather to its $+$-closure $R^+$, which has these properties by definition.
    The resulting poset will then serve as the basis for the construction of the canonical DAG $\cG_R$.
    }

\REV{\begin{definition}\label{def:R+poset}
        Let $R$ be a relation on $\pairs(X)$.
        We define the equivalence relation $\sim_{R^+}$ on $\support^+_R$ by putting, for all $p,q \in \support^+_R$,
        \[
        p\sim_{R^+} q \quad\Longleftrightarrow\quad p\ R^+\ q\ \text{ and }\ q\ R^+\ p.
        \]
        Let $[p]$ denote the equivalence class of $\sim_{R^+}$ that contains $p\in\support_R^+$, and let $Q$ denote the set of all such equivalence classes.
        Define the partial order $\le_{R^+}$ on $Q$ by putting, for all classes
        $[p]$ and $[q]$ in $Q$,
        \[
        [p]\ \le_{R^+}\ [q] \quad \Longleftrightarrow\quad p\ R^+\ q.
        \]
        We refer to the poset $(Q,\le_{R^+})$ as the \emph{quotient poset of $R^+$}.
    \end{definition}
}

\REV{The well-definedness of the equivalence relation $\sim_{R^+}$ and of the poset
$(Q,\le_{R^+})$ in Definition~\ref{def:R+poset} follows from the standard
quotient construction that turns a preorder into a partial order; see \cite[Prop.~5.2.4]{Schroder:03} 
or \cite[Lem.~5.3\&5.4]{LAVGH:25}
for a full proof.}

We now define the canonical \REVII{DAG} of a relation $R$ on $\pairs(X)$.

\begin{definition}[Canonical \REVII{DAG}]\label{def:canonG}
Let $R$ be a   relation on $\pairs(X)$ and let \REVII{$(Q,\leq_{R^+})$ be the quotient poset of $R^+$}.
The \emph{canonical \REVII{DAG}} $\cG_R$ is defined 
\REVII{as the DAG obtained from the Hasse diagram $\Hasse(Q,\leq_{R^+})$ by relabeling all vertices $[aa]\in Q$ with $a$.}
\end{definition}

\REV{Note that if a relation $R$ does not satisfy \axiom{X1}, then it may contain constraints of the form $(ab,xx)\in R$. In this case, it is possible that, for
example, $[aa]=[bb]\in Q$. Hence, in the construction of $\cG_R$, the class $[aa]$ would have to be replaced by the leaf $a$, while the same class, viewed
as $[bb]$, would have to be replaced by the leaf $b$. Thus, the construction of $\cG_R$ would no longer be well-defined.
One can circumvent this ambiguity by choosing, for each such class, a single representative leaf by which it is replaced (e.g. by imposing a linear order on $X$ and choosing the representative as the smallest element).  Nevertheless, we will only explicitly construct $\cG_R$ for relations $R$ on $\pairs(X)$ that satisfy
\axiom{X1}. Under this assumption, the ambiguity described above cannot occur.
Consequently, $\cG_R$ is a well-defined DAG and, as we shall see in Proposition~\ref{prop:CanonicalG-DAG-outsourced}, has leaf set $X$.
}

\begin{figure}
    \centering
    \includegraphics[width=0.8\textwidth]{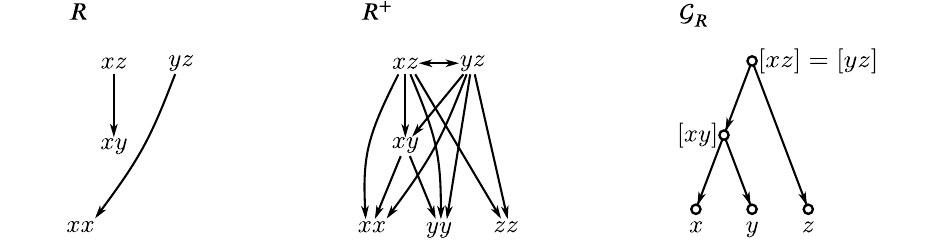}
    \caption{
    On the left we give a graphical representation of the  relation 
  $R= \{(xy,xz), (xx,yz)\}$ on \REV{$\pairs(X)$ with} $X=\{x,y,z\}$. 
  Here we 
  draw an arc $p\to q$ precisely if $q\ R\ p$. In addition, the graphical representation 
  of $R^+$ is provided where arcs $(p,p)$  are omitted for all $p\in \support_R^+$. Furthermore,
  the canonical DAG $\cG_R$ is shown. Here, $\cG_R$ realizes $R$.}
    \label{fig:working-example}
\end{figure}

An example of the canonical \REVII{DAG} $\cG_R$ of a realizable relation $R$ is provided in Figure~\ref{fig:working-example}.
We emphasize that the construction of $\cG_R$ is solely based on properties of the 
relation $R^+$. Similarly, the construction of $\cG_{R^+}$ depends only on the relation $(R^+)^+$. 
Since, by Proposition~\ref{lem:closure-axioms}, $(R^+)^+ = R^+$, it follows that the construction of $\cG_{R^+}$ is determined by the relation $R^+$, just as for $\cG_R$.
Hence, we obtain
\begin{observation}\label{obs:G_R=G_R+}
    For all relations $R$ we have $\cG_R=\cG_{R^+}$.
\end{observation}

As we shall see in Theorem~\ref{thm:char}, if $R$ can be realized by some DAG, then it is realized by $\cG_R$.  To this end, we first give two useful properties of $\cG_R$ in the case that $R$ 
satisfies Condition \axiom{X1}.

\begin{proposition}\label{prop:CanonicalG-DAG-outsourced}
Let $R$ be a relation on $\pairs(X)$ that satisfies Condition \axiom{X1}.
Then, the canonical \REVII{DAG} $\cG_R$  is a 
DAG on $X$ such that $[ab]=\lca_{\cG_R}(ab)$ for all vertices $[ab]$ of $\cG_R$ with $a\neq b$
and $x=\lca_{\cG_R}(xx)$ for all $x\in X$. In particular, $\cG_R$ is 2-lca-relevant, phylogenetic
and realizes $R^+$.
\end{proposition}
\begin{proof} 
\REVII{Let $R$ be a relation on $\pairs(X)$ that satisfies Condition \axiom{X1} an let $(Q,\le_{R^+})$ be the quotient poset of $R^+$.
We start by showing that the leaf set of $G\coloneqq \cG_R$ is $X$. 
Since $R$ satisfies \axiom{X1}, Lemma~\ref{lemma:X1-R-R+} ensures that $R^+$ satisfies \axiom{X1}. Hence,} no class $[xx]\in Q$ 
can coincide with a class $[ab]\in Q$ for which $ab\neq xx$. 
By definition, $R^+$ is $\support_R^+$-reflexive and, thus, $[xx]\in Q$ for all $x\in X$. 
Hence, the two sets $Q_1=\{[xx]\colon x\in X\}$ and $Q_2 =\{[ab]\in Q\colon a\neq b\}$
form a bipartition of $Q$, i.e., $Q =Q_1\cup Q_2$ and $Q_1\cap Q_2 =\emptyset$. 
By the latter arguments, we have $x,y\in X$ and $x\neq y$ if and only if  $[xx]\neq [yy]$ and $[xx],[yy]\in Q_1\subseteq Q$. 
Thus, there is a 1-to-1 correspondence between the vertices in $X$ and the classes in $Q_1$. 
Moreover, \REVII{since $R^+$ satisfies} \axiom{X1},  $(ab, xx)\not\in R^+$ for all $a,b,x\in X$ with $ab \neq xx $.
Hence, $[ab] \not \le_{R^+} [xx]$ for all $a,b,x\in X$ with $ab \neq xx$. 
Therefore, arcs of the form $x\to v$ for any $x\in X$ \REVII{do not exist in} $G$. Moreover, 
by construction, $G$ does not contain arcs $(v,v)$ for any $v\in V(G)$. In summary,
the set $X$ forms a subset of leaves in $G$. 
\REVII{It remains to show that all vertices $[xy]\in Q$ with 
$x\neq y$ are non-leaf vertices in $G$. To this end, observe
that cross-consistency of $R^+$ ensures that
$[xx]\leq_{R^+} [xy]$ for all classes $[xy]$ in $Q$, cf.\ Observation~\ref{obs:always_xx<xy}. 
Hence, for all such classes $[xy]$ with $x\neq y$, 
the vertex $x$ satisfies $x\prec_{G} [xy]$. 
Thus, classes $[xy]$ of $Q$ with $x\neq y$
must be non-leaf vertices in $G$. }
In summary, $X$ is precisely the leaf set  of $G$.


Now let $v\in V(G)$. Suppose that $v=[ab]$ for some class $[ab]\in Q$ with $a\neq b$. 
We show that $v$ is the unique least common ancestor of $a$ and $b$ in $G$.
\REVII{By the previous arguments,}
 $[ab]$ is a common ancestor of $a$ and $b$.
Now suppose that $w \in V(G)$ is a common ancestor of $a$ and $b$
with $w \neq [ab]$. Then $w$ must be in $Q_2$, i.e., $w=[xy]\in Q$ with $x\neq y$. 
Since $w$ is a common ancestor of $a$ and $b$, there are directed paths $[xy]\leadsto a$ and $[xy]\leadsto b$ in $G$.  By construction of $G$ and since $\le_{R^+}$ is transitive (cf.\ Definition~\ref{def:R+poset}), 
we must have $[aa]\le_{R^+} [xy]$
	and $[bb]\le_{R^+} [xy]$. By definition of $\le_{R^+}$, we have $aa\ R^+\ xy$ and $bb\ R^+\ xy$. Since, $[ab]\in Q$ we have, by construction of $\sim_{R^+}$, $ab\in\support_R^+$. 
    By Theorem~\ref{thm:rules}, $\support_R^+ =\support_{R^+}$ and, therefore, $ab\in\support_{R^+}$. This together with that fact that $R^+$ is cross-consistent
    implies that $ab\ R^+\ xy$. 
	Again, by definition of $\le_{R^+}$, we have $[ab]\le_{R^+} [xy] =w$. \REVII{By definition of $G$, we have therefore} 
	$[ab]\preceq_G w$.
Since $w=[xy]$ was an arbitrary common ancestor of $a$ and $b$,
$[ab] \preceq_G w$ holds for \emph{all} 
common ancestors of $a$ and $b$. Hence $[ab]$ is the \emph{unique} least common ancestor of $a$ and $b$.
Moreover, if $v=a$ for some $a\in X$,  the unique least common ancestor $\lca_G(aa)$ is $a$ since $a\in X$ is a leaf. 
In summary, $G$ is a DAG in which $\lca_G(aa)=a$ for all $a\in X$ and 
$\lca_G(ab)=[ab]$ for all $[ab]\in Q$ with $a\neq b$. 

By the latter arguments and since each vertex in $G$ is either of the form $[ab]\in Q_2$ or $x$ for some $[xx]\in Q_1$, 
it follows that, for each vertex $v$ in $G$, there are $x,y\in X$ such the $v=\lca_G(xy)$. 
Therefore, $G$ is 2-lca-relevant and, in particular, lca-relevant. 
This observation allows us to apply \cite[L.~3.10]{HL:24} to conclude that $G$ is phylogenetic.

We show now that $G$ realizes $R^+$. By definition, $R^+$ is transitive and thus, 
$\tc(R^+)=R^+$. Let $(ab,cd)\in R^+ = \tc(R^+)$. Suppose first that $(cd,ab)\notin \tc(R^+)=R^+$. 
Hence, $[ab]\neq [cd]$ and $[ab]\le_{R^+} [cd]$. \REVII{If $a\neq b$, the vertex $v = [ab]$ exists in $G$ and if $a=b$, the vertex $v=a\in X$
exists in $G$. Either way, $[ab]\neq [cd]$, $[ab]\le_{R^+} [cd]$ and the construction of $G$ ensures that $\lca_G(ab)=v\prec_G[cd]=\lca_G(cd)$. In other words,} \axiom{I1} holds for $G$ and $R^+$. 
If $(cd,ab)\in \tc(R^+)=R^+$, then $ab\sim_{R^+} cd$ and thus, $[ab]=[cd]$. 
This together with the latter arguments implies that $\lca_G(ab) = \lca_G(cd)$ holds. 
Thus, \axiom{I2} holds for $G$ and $R^+$. In summary, $G$ realizes $R^+$.
\end{proof}

\REV{Before proceeding, we note that the part of the proof of Proposition~\ref{prop:CanonicalG-DAG-outsourced}, which shows that the canonical construction in Definition~\ref{def:R+poset} yields a DAG on $X$ whenever $R$ satisfies \axiom{X1}, would still work if $R^+$ was replaced by some other reflexive and transitive relation containing $R$. However, the $+$-closure incorporates cross-consistency,
which is crucial to show that $\cG_R$ satisfies $\lca_G(ab)=[ab]$ for all vertices $ab$ of $\cG_R$.
We illustrate this on the following example.}

\begin{figure}
    \centering
    \includegraphics[width=0.8\textwidth]{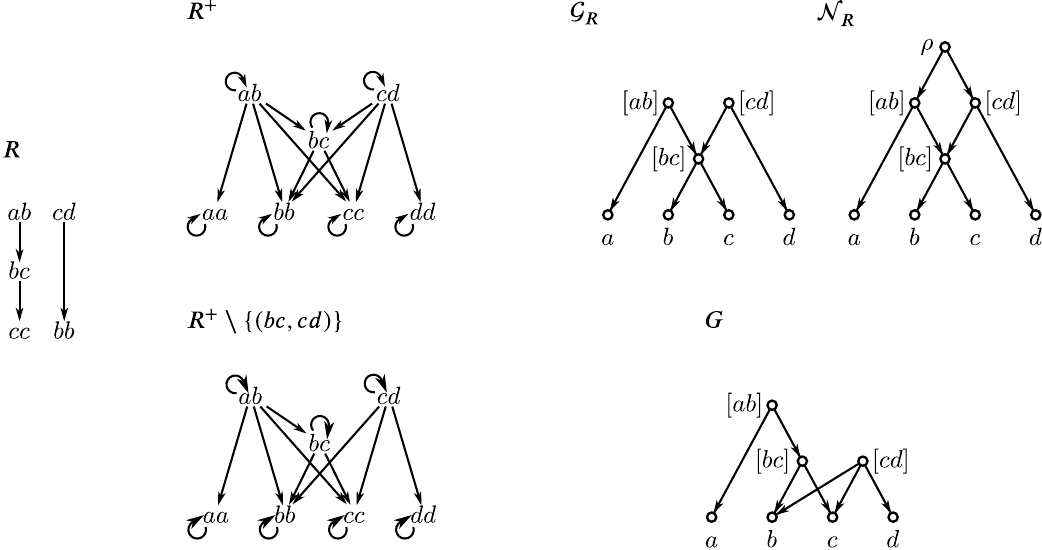}
    \caption{\REV{On the left we give a graphical representation of the relation $R=\{(cc,bc),(bc,ab),(bb,cd)\}$, where we draw an arc $p\to q$ precisely if $q\ R\ p$. Moreover, shown are the two relations $R^+$ (c.f.\ Figure~\ref{fig:stepwise-closure}) and $R'\coloneqq R^+\setminus\{(bc,cd)\}$ as well as three different DAGs related to these two relations; see Example~\ref{ex:cross-consistency} and Section~\ref{sec:canon-N} for details.}}
    \label{fig:cross-consistency}
\end{figure}

\REV{
\begin{example}[Cross-consistency]\label{ex:cross-consistency}
    Consider the relation $R=\{(cc,bc),(bc,ab),(bb,cd)\}$ in Figure~\ref{fig:cross-consistency} which satisfies \axiom{X1},
    see also Figure~\ref{fig:stepwise-closure}. 
    The corresponding canonical DAG $\cG_R$ satisfies, as ensured by Proposition~\ref{prop:CanonicalG-DAG-outsourced}, that $\lca_{\cG_R}(xy)=[xy]$ for each $xy\in\{ab,bc,cd\}$. In fact, it is easily verified that $\cG_R$ indeed realizes $R$. Figure~\ref{fig:cross-consistency} also shows a graphical representation of $R'\coloneqq R^+\setminus\{(bc,cd)\}$, which is a $\support_R^+$-reflexive and transitive relation on $\pairs(X)$. Note that $R'$ is not cross-consistent, since $(bb,cd),(cc,cd)\in R'$, $bc\in\support_{R'}^+$ and $(bc,cd)\notin R'$. However, since $R'$ is transitive and $\support_R^+$-reflexive, one may mimic Definitions~\ref{def:R+poset} and \ref{def:canonG} to construct the DAG $G$ on $X$ shown at the bottom left of Figure~\ref{fig:cross-consistency}. Note, however, that in $G$ we have $\LCA_G(bc)=\{[bc],[cd]\}$, so $\lca_G(bc)$ is undefined and $G$ does not realize $R$.
\end{example}
}

\REV{Example~\ref{ex:cross-consistency} shows that the cross-consistency of $R^+$ is essential: it ensures that the
poset constructed in Definition~\ref{def:R+poset} has the additional  structural
properties required to build a DAG $G$ on $X$ that realizes $R$, provided that
$R$ is realizable.
}\REV{Moreover,} we note that Proposition~\ref{prop:CanonicalG-DAG-outsourced} ensures that the canonical DAG $\cG_R$  realizes $R^+$ if $R$ satisfies \axiom{X1}. 
This, however, does not imply that $R$ is realizable, as shown in the following example. 

\begin{example}[Relations $R$ that are not realizable, although $R^+$ is realizable]\label{exmpl:RnotR+is}
Let $R =\{(xx,ab), (yy,ab), (ab,xy)\}$ be the relation on \REV{$\pairs(X)$ with} $X = \{x,y,a,b\}$
as provided in Example~\ref{exmpl:not-subset-real}. As argued in   
Example~\ref{exmpl:not-subset-real}, $R$ is not realizable by any DAG $G$.
However, $R$ clearly satisfies
\axiom{X1}, which together with Proposition~\ref{prop:CanonicalG-DAG-outsourced} 
implies that $\cG_R$ realizes $R^+$. Hence, $R^+$ is realizable, although $R$ is not. 
\end{example}

 Example~\ref{exmpl:RnotR+is} shows that \axiom{X1} 
is not sufficient on its own for a relation $R$ to be realizable.
In particular, contraposition of Lemma~\ref{lem:no-collapse}
shows that also \axiom{X2} alone is insufficient to characterize 
realizable relations. We now prove the main result of this section.

\begin{theorem}\label{thm:char}
For a relation $R$ on $\pairs(X)$ the following statements are equivalent:
\begin{enumerate}
  \item $R$ is realizable.
  \item $R$ satisfies \axiom{X1} and \axiom{X2}.
  \item $R$ is realized by its canonical DAG $\cG_R$.
\end{enumerate}
\end{theorem}
\begin{proof}
If Condition (1) is satisfied, then  Lemma~\ref{lem:no-collapse} implies that 
 $R$ satisfies \axiom{X1} and \axiom{X2}, i.e., Condition (2) holds. 
 
Assume now that Condition (2) is satisfied. 
Consider the canonical DAG $G\coloneqq \cG_R$ of $R$. 
Proposition~\ref{prop:CanonicalG-DAG-outsourced}  implies
that $G$ is DAG on $X$ in which $\lca_G(aa)=a$ for all $a\in X$ and 
$\lca_G(ab)=[ab]$ for all $[ab]\in Q$ with $a\neq b$.
We show that $R$ is realized by $G$. Let $ab,xy\in \support_R^+$ and suppose that 
\REVII{$ab\ R\ xy$. Since $R\subseteq R^+$, we indeed also have $ab\ R^+\ xy$.}

Assume, first that $(xy,ab)\not\in \tc(R)$. 
By \axiom{X2}, $(xy,ab)\notin R^+$. 
Therefore, $[ab] \neq [xy]$ and  $[ab] \leq_{R^+} [xy]$.
If $a\neq b$, the vertex $v = [ab]$ exists in $G$ and if $a=b$, the vertex $v=a\in X$
exists in $G$. In either case, 
\REVII{$v\prec_G[xy]$} and, by Proposition~\ref{prop:CanonicalG-DAG-outsourced},
$\lca_G(ab)=v$ holds. Moreover, Proposition~\ref{prop:CanonicalG-DAG-outsourced} implies $ [xy]=\lca_G(xy)$.
In summary, $\lca_G(ab)=v \prec_G\ [xy]=\lca_G(xy)$. Hence, Condition \axiom{I1} holds.

Assume, now that $(xy,ab)\in \tc(R)$. Recall that $\tc(R)\subseteq R^+$, so
we have $(ab,xy),(xy,ab)\in R^+$. Hence $ab$ and $xy$ are in the same $\sim_{R^+}$-class, 
i.e., $[ab]=[xy]$ refer to the same vertex in $G$. 
By Proposition~\ref{prop:CanonicalG-DAG-outsourced}, we have 
$[xy]=\lca_G(xy)$ and $[ab]=\lca_G(ab)$ which together with $[ab]=[xy]$ implies that 
$\lca_G(xy)=\lca_G(ab)$. Thus, Condition \axiom{I2} holds. In summary, $G$ realizes $R$.
Hence, Condition (2) implies Condition (3).

Clearly, Condition (3) implies (1). In summary, the three statement (1), (2) and (3)
are equivalent.
\end{proof}

\begin{corollary}\label{cor:R=R+=>(X2)}
    Let $R$ be a relation such that $R = R^+$. Then $R$ satisfies \axiom{X2}. 
        In particular, $R$ is realizable if and only if it satisfies \axiom{X1}. 
\end{corollary}
\begin{proof}
        Let $R$ be a relation such that $R = R^+$.
        Since $R^+$ is transitive and $R = R^+$, we have $\tc(R) = \tc(R^+) = R^+$. 
        Therefore, Condition~\axiom{X2} is trivially satisfied. 
        Thus, if $R$ satisfies \axiom{X1}, Theorem~\ref{thm:char} implies that $R$ is realizable.
        Conversely, if $R$ is realizable, then $R$ satisfies \axiom{X1} by Theorem~\ref{thm:char}.
\end{proof}

The latter results also allow us to characterize strictly realizable relations as follows.

\begin{theorem}\label{thm:char-strict}
For a relation $R$  the following statements are equivalent:
\begin{enumerate}
  \item $R$ is strictly realizable.
  \item $R$ satisfies \axiom{X1} and  \axiom{X2} and  $\tc(R)$ is asymmetric.
  \item $R$  is realized by its canonical DAG $\cG_R$ and $\tc(R)$ is asymmetric.
  \item $R\subseteq\srel_{\cG_R}$ for the canonical DAG $\cG_R$ of $R$.
  \item $R\subseteq\srel_G$ for some DAG $G$.
\end{enumerate}
\end{theorem}
\begin{proof}
    By Theorem~\ref{thm:char} together with Lemma~\ref{lem:str-real=>real}, Condition~(1), (2) and (3)
    are equivalent. Suppose that Condition~(3) holds. In this case,  Lemma~\ref{lem:str-real=>real} implies that 
    $R$ is strictly realized by $\cG_R$. By Lemma~\ref{lem:Rsubset-srelG}, $R\subseteq\srel_{\cG_R}$ and thus, 
    Condition~(4) is satisfied.
    Clearly, Condition (4) implies (5). By Lemma~\ref{lem:Rsubset-srelG}, Condition (5) implies Condition (1).
\end{proof}

\begin{figure}
  \centering
  \includegraphics[width=0.8\textwidth]{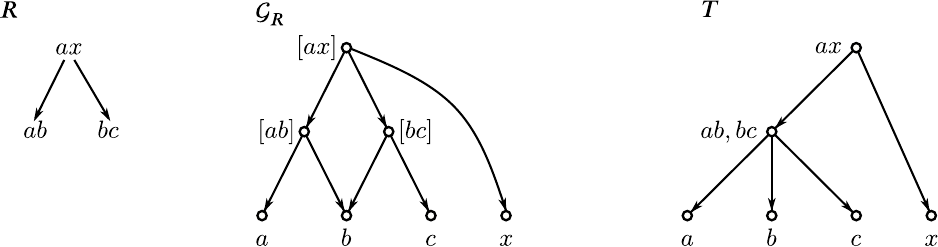}
  \caption{On the left we give a graphical representation of the relation relation $R = \{(ab,ax),(bc,ax)\}$. 
  Here we draw an arc $p\to q$ precisely if $q\ R\ p$. In addition, the canonical DAG $\cG_R$ and a tree $T$ that both realize $R$ are shown.     
    Both $\cG_R$  and $T$ are minimal, however, only $T$ is minimum. 
  Each pair $p \in \{ab, bc,ax\}$ in the drawing of $T$ marks the corresponding vertex $\lca_{T}(p)$.
  }
  \label{fig:non-min}
\end{figure}

We conclude this section with a brief discussion of ``minimal'' and ``minimum'' DAGs
realizing a relation $R$. Observe first that the canonical DAG $\cG_R$ of a realizable relation $R$ is, in general, not the only DAG that realizes $R$. Moreover, the Examples in Figure~\ref{fig:non-min} and Figure~\ref{fig:exmpl-tree-reg} show that $\cG_R$ is not necessarily \emph{minimum}, i.e., it does not  always have the minimum number of vertices among all DAGs $G$ that realize $R$. 

Although the canonical DAG is not necessarily minimum with respect to the number of vertices, 
it may still be minimal in a structural sense. 
To formalize this, we call a DAG $G$ that realizes a relation $R$  \emph{minimal} if every DAG $H$ obtained from $G$ by a sequence of arc contractions\footnote{In an arc contraction, an arc $(u,v)$ 
    is removed from $G$
    and $u$ and $v$ are identified, and then any resulting multi-arcs are 
    suppressed to give a single arc.} fails to realize $R$.
However, as the following example illustrates, $\cG_R$ is not always minimal.

\begin{figure}
  \centering
  \includegraphics[width=0.8\textwidth]{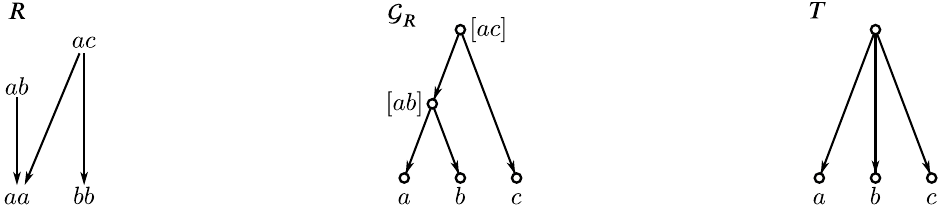}
  \caption{On the left we give a graphical representation of the relation relation 
  $R= \{(aa,ac), (bb,ac),(aa,ab)\}$ on $\pairs(X)$ with $X=\{a,b,c\}$ as in Example~\ref{exmpl-not-minimal}. 
  Here we 
  draw an arc $p\to q$ precisely if $q\ R\ p$. In addition, the canonical DAG $\cG_R$ a DAG $T$ is shown. Both $T$ and $\cG_R$ are trees that realize $R$. 
  Moreover, contracting the arc $[ac]\to [ab]$ in $\cG_R$ yields $T$. Hence, $\cG_R$ is not a ``minimal'' DAG realizing $R$. }
  \label{fig:non-minimum}
\end{figure}

    \begin{example}[The canonical DAG is not minimal]\label{exmpl-not-minimal}
        Let $R= \{(aa,ac), (bb,ac),(aa,ab)\}$ be a relation on $\pairs(X)$, where $X=\{a,b,c\}$. 
        One easily verifies that the star-tree $T$ on $X$ realizes $R$. 
        By cross-consistency and since $(aa,ac), (bb,ac)\in R$ and $ab\in \support_{R^+}$
        it follows that $(ab,ac)\in R^+$. It is also easy to see that $(ac,ab)\notin R^+$.
        The canonical DAG $\cG_R$ is shown in Figure~\ref{fig:non-minimum}. 
        One easily verifies that $T$ can be obtained from $\cG_R$ by contracting 
        the arc $ac\to ab$.
        Hence, \REVII{$\cG_R$ is not minimal.}
    \end{example}

\section{The Canonical Network:  Properties and Algorithmic Construction}
\label{sec:canon-N}

In the previous section, we saw how to realize a relation $R$ by its canonical DAG $\cG_R$. 
In this section, we focus on realizing a relation via a canonical \emph{network} $\cN_R$. 
This network is  obtained by a simple modification of $\cG_R$. 
Interestingly, we will show that \REVII{$\cG_R$ and }$\cN_R$ are regular, meaning that they are closely related 
to the Hasse diagram of their underlying set system, which arises by taking the 
descendants in $X$ of vertices in \REVII{$\cG_R$, respectively,} $\cN_R$. 
In addition, we present a polynomial-time algorithm for computing both $\cG_R$ and $\cN_R$.

We start by relating DAGs to set systems.
 For a DAG $G$ on $X$ and a vertex $v\in V(G)$ we put $\CC_G(v)=\{x\in X\mid x\preceq_G v\}$, i.e. the set of leaves that are descendants of $v$.
This set is called the \emph{cluster} of $v$, and we let $\mathfrak{C}_G=\{\CC_G(v)\mid v\in V(G)\}$ denote the set system on $X$ comprising all clusters of $G$. Clearly $\mathfrak{C}_G$ is grounded for all DAGs $G$. In the opposite direction, we can also associate 
\REVII{the Hasse diagram $\Hasse(\mathfrak{C},\subseteq)$
to the poset $(\mathfrak{C},\subseteq)$ where 
$\mathfrak{C}\subseteq \mathcal{P}(X)$ is a set system
and $\subseteq$ denotes the usual subset-relation. 
Hence, $\Hasse(\mathfrak{C},\subseteq)$ is a DAG with}
vertex set $\mathfrak{C}$ and an arc from $A\in\mathfrak{C}$ to $B\in\mathfrak{C}$
if (i) $B\subsetneq A$ and (ii) there is no $C\in\mathfrak{C}$ with $B\subsetneq C\subsetneq A$. We
note that $\Hasse(\mathfrak{C},\subseteq)$ is also known as the \emph{cover digraph} of $\mathfrak{C}$
\cite{Baroni:05}. We now give a key definition.

\begin{definition}[{\cite{Baroni:05}}]
  \label{def:regular-N}
  A DAG $G=(V,E)$ is \emph{regular} if the map
  $\varphi\colon V\to V(\Hasse)$ defined by  $v\mapsto \CC_G(v)$ is an
  isomorphism between $G$ and $\Hasse \REVII{\coloneqq \Hasse(\mathfrak{C},\subseteq)}$.
\end{definition}

We now define the canonical network for a relation.

\begin{definition}[Canonical Network]\label{def:canon-network}
Suppose that $R$ is a realizable relation on $\pairs(X)$. Let $\cG_R$ be the canonical DAG which, by Theorem~\ref{thm:char} realizes $R$. 
We construct the \emph{canonical network $\cN_R$} from $\cG_R$ as follows. 
\REVII{\begin{itemize}
	\item If $\cG_R$ is a network, then put $\cN_R\coloneqq \cG_R$.
    
	\item Else, $\cG_R$ has roots $\rho_1,\dots, \rho_k$, $k\geq 2$, and we obtain $\cN_R$
			by first adding a new vertex $\rho$ and then adding the arc  $\rho\to \rho_i$,  
            for each $i\in \{1,\dots,k\}$.
\end{itemize}}
\end{definition}

\begin{figure}
    \centering
    \includegraphics[width=0.8\textwidth]{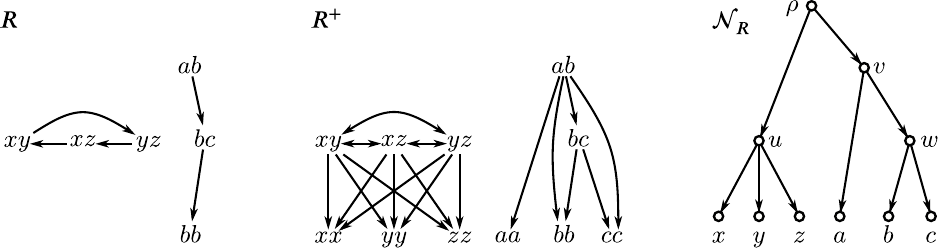}
    \caption{On the left we give a graphical representation of the relation relation $R$.
                 Here we draw an arc $p\to q$ precisely if $q\ R\ p$. In addition, the graphical representation of $R^+$ is provided where arcs $(p,p)$  are omitted for all $p\in \support_R^+$. Furthermore, the canonical network $\cN_R$ is shown. In this example $\cG_R$, which is distinct from $\cN_R$, is obtained from $\cN_R$ by removal of the root $\rho$ and its incident arcs.
                 In $\cN_R$, we have $u=[xy]=[xz]=[yz]$, $v=[ab]$, and $w=[bc]$.
    } 
    \label{fig:working-example-2}
\end{figure}

\begin{figure}
  \centering
  \includegraphics[width=0.8\textwidth]{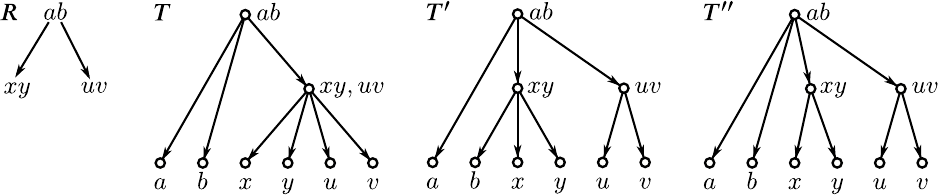}
  \caption{Three non-isomorphic phylogenetic trees $T$, $T'$ and $T''$ on $X=\{a,b,x,y,u,v\}$ that all  realize
            the relation $R$ with $xy\ R\ ab$ and $uv \ R\ ab$.  
            Each pair $p \in \{ab, xy, uv\}$ in the drawings marks the corresponding vertex $\lca_{T^*}(p)$ in the tree $T^* \in \{T, T', T''\}$.
            Here, $T''\simeq \cG_R=\cN_R$. Moreover, all three trees are regular
            and neither of them can be obtained from the other by contraction of arcs, i.e., 
            they are minimal. However, only $T$ has a minimum number of vertices among all DAGs realizing $R$. 
            }
  \label{fig:exmpl-tree-reg}
\end{figure}

\REV{Examples of canonical networks that differ from the underlying canonical DAG are shown in Figure~\ref{fig:cross-consistency} and Figure~\ref{fig:working-example-2}.}  We now show that,
as well as some enjoying other properties, the canonical network is always regular.
Note that, even so, regular DAGs or networks that realize $R$ are not necessarily uniquely determined -- 
see e.g. Figure~\ref{fig:exmpl-tree-reg}.

\begin{proposition}\label{prop:srel-realized-by-N}
If $R$ is  realizable, then \REVII{$\cG_R$ is regular and}
the canonical network $\cN_R$ realizes $R$. Moreover, 
 $\cN_R$ is regular and phylogenetic and satisfies $[ab]=\lca_{\cN_R}(ab)$ for all vertices $ab\in \support_{R}^+$ with $a\neq b$ and $x=\lca_{\cN_R}(xx)$ for all $x\in X$.
\end{proposition}
\begin{proof}
    Suppose that $R$ is a realizable relation on $\pairs(X)$.
    In the following, we put $G\coloneqq \cG_R$ and $N\coloneqq \cN_R$.
   	\REVII{By Proposition~\ref{prop:CanonicalG-DAG-outsourced}, $G$ is 2-lca-relevant and, hence,  lca-relevant.
   	Moreover, by construction, $G$ does not contain shortcuts.
   	The latter arguments together with  \cite[Thm.~4.10]{HL:24} 
   	imply that $G$ is regular. This and Proposition~\ref{prop:CanonicalG-DAG-outsourced}
   	implies that $N=G$ satisfies the stated conditions whenever $G$
   	is already network. 
   	}

Suppose now that \REVII{$G$} is not a network. In this case, $G$ has roots 
$\rho_1,\dots, \rho_k$ with $k\geq 2$. We first argue that $\CC_{G}(\rho_i)\neq X$, $1\leq i \leq k$. 
Assume, for contradiction, that $\CC_{G}(\rho_i)= X$ for some $i$. In this case, $\CC_{G}(v)\subseteq 
\CC_{G}(\rho_i)$  for all $v\in V(G)$. Since $G$ is regular,  \cite[L.~4.7]{HL:24} implies that 
$G$ satisfies the so-called path-cluster-comparability (PCC) property which, in particular, implies
that $\rho_i$ and $\rho_j$ must be $\preceq_{G}$-comparable for all $j$; a contradiction
to $\rho_i$ and $\rho_j$ being roots for at least two distinct $i$ and $j$. Consequently, $\CC_{G}(\rho_i)\neq X$ for each $1\leq i \leq k$. 

Now, recall that we construct $N$ by taking $G$ and 
adding a new vertex $\rho$ and arcs $\rho\to \rho_i$ for all $1\leq i \leq k$. It is easy to verify that 
$N$ is a phylogenetic network. To show that $N$ is regular, we show first that $N$ is lca-relevant. 
By Proposition~\ref{prop:CanonicalG-DAG-outsourced}, $G$ is 2-lca-relevant, i.e., for all $v\in V(G)$
there are leaves $x,y\in X$ (not necessarily distinct) such that $v = \lca_{G}(xy)$. It is straightforward to verify that
$v = \lca_{N}(xy)$ whenever $v\in V(G) = V(N)\setminus \{\rho\}$ and $v = \lca_{G}(xy)$.
We show now that $\rho=\lca_N(X)$.
As argued above, all children of $\rho_i$ of $\rho$ in $N$, 
satisfy $\CC_N(\rho_i)\neq X$.
 Moreover, any descendant $v$ of $\rho_i$ must satisfy $\CC_N(v)\subseteq \CC_N(\rho_i)$ (cf. e.g. \cite[L.17]{Hellmuth2023}). 
 Hence, $\CC_N(v)\neq X$ for all $v\in V(N)\setminus\{\rho\}$. 
Since $\rho$ is an ancestor of all leaves in $X$ and since there is no
vertex $v\prec_N  \rho$ with  $\CC_N(v) = X$, it follows that $\lca_N(X)=\rho$. 
In summary, for all vertices $v\in V(N)$  there is a subset $A\subseteq X$
such that $v=\lca_N(A)$. Therefore,
$N$ is lca-relevant. 

We show now that $N$ is shortcut-free. Observe first that 
the subgraph $G$ of $N$ is shortcut-free. Hence any shortcut in $N$
must be of the form $\rho\to v$ for some $v\in V(N)$. Since $\rho$ is only adjacent to $\{\rho_1,\dots,\rho_k\}$, 
a shortcut must be of the form $\rho\to \rho_i$ for some $1 \le i \le k$. However,  $\rho\to \rho_i$ being a shortcut
requires an alternative $\rho \rho_i$-path not containing the arc $\rho\to \rho_i$. Since, by construction, 
no such path exists in $N$, $N$ is shortcut-free. Since  $N$ is, in addition, lca-relevant, 
we can apply \cite[Thm.~4.10]{HL:24} to conclude that $N$ is regular.

Finally, recall that $\lca_{G}(ab)$ is well-defined for all $ab\in \support^+_R$.
Note that, if $ab\in \support^+_R$, then $\lca_{G}(ab)$ is  contained in $V(G) = V(N)\setminus \{\rho\}$.
Moreover, $G$ is an induced subgraph of $N$. It is now straightforward to verify that $N$ realizes $R$
and that $[ab]=\lca_{N}(ab)$ for all vertices $ab\in \support_R^+$ with $a\neq b$ and $x=\lca_{N}(xx)$ for all $x\in X$. 
\end{proof}

We conclude this section by presenting the polynomial-time algorithm
which checks whether or not a relation $R$ is realizable and, if so, returns 
the canonical DAG and network realizing $R$. 
The algorithm is presented in Algorithm~\ref{alg:REAL} and its correctness
and runtime analysis is given in the following result.

\begin{algorithm}
  \caption{\textsc{Realizability of a Relation $R$}}
  \label{alg:REAL}
  \begin{algorithmic}[1]
  \Require  A binary relation  $R$ on $\pairs(X)$
    \Ensure   Returns the canonical DAG $\cG_R$ and network $\cN_R$ realizing $R$ if $R$ is realizable and, otherwise, \texttt{false} is returned.
    \State Compute $\support_{R}^+$ \label{l:supp}
    \State Compute $R^+$ \REV{with Algorithm~\ref{alg:Rplus}.} \label{l:R+}
	 \If{$R$ satisfies Conditions \axiom{X1} and \axiom{X2}}   \label{l:X12}
     \State \REVII{Compute the poset $(Q,\le_{R^+})$ as in Definition~\ref{def:R+poset} \label{l:poset}}
	 \State Compute the canonical DAG $\cG_R$ according to Definition~\ref{def:canonG} \label{l:GR}
	 \State Compute the canonical network $\cN_R$ according to Definition~\ref{def:canon-network} \label{l:NR}
 	 \State \Return $\cG_R$ and $\cN_R$
 	 \Else \ 	\Return \texttt{false} \EndIf
  \end{algorithmic}
\end{algorithm}

\begin{theorem}\label{thm:alg:char}
For a relation $R$ on $\pairs(X)$, verifying whether $R$ is realizable and, if so, constructing a DAG on $X$ that realizes $R$ can be done in polynomial time in $|X|$ using Algorithm~\ref{alg:REAL}. In particular, \REVII{both $\cG_R$ and $\cN_R$,} can be constructed in polynomial time,
where $\cG_R$ is the canonical DAG and $\cN_R$ the canonical network of $R$. 
\end{theorem}
\begin{proof}
Let $R$ be a relation on $\pairs(X)$. \REV{By Theorem~\ref{thm:R+-polytime}, 
Algorithm~\ref{alg:Rplus} correctly computes $R^+$.}
By Theorem~\ref{thm:char}, 
$R$ is realizable if and only if it satisfies \axiom{X1} and \axiom{X2}. 
It is straightforward to verify that \texttt{false} is correctly returned by Algorithm~\ref{alg:REAL} 
whenever $R$ is not realizable by any DAG~$G$. 
Lines~\ref{l:X12} to Line \ref{l:GR} mimic the steps used to construct the canonical DAG~$\cG_R$  (see also the corresponding lemmas and definitions referenced in Algorithm~\ref{alg:REAL}). 
Since $\cG_R$ is computed only when $R$ satisfies \axiom{X1} and \axiom{X2}, 
Theorem~\ref{thm:char} ensures that $\cG_R$ is indeed a DAG that realizes~$R$, whenever $R$ is realizable.

For the runtime, observe first that $\support_{R}^+$ can clearly be 
determined in polynomial time  in $|X|$  (Line \ref{l:supp}). 
Moreover, \REV{Theorem~\ref{thm:R+-polytime}} implies that $R^+$ can also be determined
in polynomial time in $|X|$ (Line~\ref{l:R+}). As outlined at the end the proof of 
\REV{Theorem~\ref{thm:R+-polytime}}, we have $|R| \in O(|X|^4)$. By similar arguments, $|R^+| \in O(|X|^4)$.
We can thus assume that $R$ is provided as an $|X|^2 \times |X|^2$ matrix
$A$ with entries $0$ and $1$ such that $A_{uv,xy}=1$ if and only of $uv\ R\ xy$. 
Hence, $A$ is the adjacency matrix of a directed graph, here called $G(A)$. 
\REV{A similar matrix representation $A^+$ of $R^+$ can also be assumed.} 
To verify $\axiom{X1}$ in Line~\ref{l:X12}, we have to check 
for all $x\in X$ and the $O(|X|^2)$ sets $ab\in \pairs(X)$ with $a\neq b$ if \REVII{$(ab,xx)\in R$} or not. 
Since we can check in constant time as whether or not  \REVII{$(ab,xx)\in R$}  based on the matrix $A$, 
 $\axiom{X1}$ can be verified in  $O(|X|^3)$ time. 
 To check $\axiom{X2}$ in Line~\ref{l:X12} we must verify if, for all  \REVII{$(ab,xy)\in R$} and  $(xy, ab)\not\in \tc(R)$,
 we have $(xy,ab)\notin R^+$.  We can compute the transitive closure of $R$ in the graph $G(A)$ with the Floyd-Warshall Algorithm
 in polynomial time in $|X|$ \REV{ and check for membership in $R^+$ with the matrix $A^+$}. Hence, $\axiom{X2}$  can be verified in polynomial time.  
 Computing \REVII{the poset $(Q,\le_{R^+})$ in Line~\ref{l:poset}}  can be achieved in 
 polynomial time in $|X|$, since its computation requires only a  polynomial number of 
 comparison of elements in $R^+$.
Since $\cG_R$, by Proposition~\ref{prop:CanonicalG-DAG-outsourced}, is 2-lca-relevant, 
it follows that each vertex in $\cG_R$ is the unique least common ancestor of at least two 
vertices out of $|\pairs(X)|$ possible leaf combinations. We can, thus, conclude that 
the number of vertices and arcs in $\cG_R$ is bounded by $|\pairs(X)|\in O(|X|^2)$ and  $O(|X|^4)$, respectively.
\REVII{With this it is now easy to verify that the graph
 $\cG_R$ (Line~\ref{l:GR}) as specified in Def.~\ref{def:canonG} 
 as well as $\cN_R$ (Line~\ref{l:NR})
 can be computed in polynomial time in $|X|$.}
 Hence, Algorithm~\ref{alg:REAL} runs in polynomial time in $|X|$. 
\end{proof}

It is easy to check that verifying if a relation is asymmetric can be performed in 
polynomial time. This together with Theorem~\ref{thm:char-strict} and Theorem~\ref{thm:alg:char} implies the following result. 

\begin{theorem}
Verifying whether a relation $R$ is strictly realizable and, if so, constructing a DAG that strictly realizes $R$ can be done in polynomial time in $|X|$.    
\end{theorem}

\section{The Classical Closure and Closed Relations}
\label{sec:proof-closure}

In previous sections, we defined the closure $R^+$ of a relation $R$ and
then showed how to use this to construct the canonical DAG for a realizable relation. However, a more natural way to define the closure for a realizable relation is to 
follow the approach used for so-called triplets or quartets in 
trees \cite{SH:18,MaayanLevy24,Bryant97,BS:95} \REV{(for further discussion see Section~\ref{sec:outlook})}.
More specifically, we define the closure $\cl(R)$ as the intersection of all relations $\rel_G$ 
over DAGs $G$ that realize $R$. 
In this section we will show that, somewhat surprisingly, these two ways of taking 
closures are actually equivalent (see Theorem~\ref{thm:closure-classic}).

We begin by formally defining the alternative notion of closure.

\begin{definition}[Classical closure]\label{def:closure}
    Let $R$ be a   relation on $\pairs(X)$ that is realizable
    and let $\mathfrak{G}$ denote the set of all DAGs $G$ on $X$ that realize $R$. 
    Then, we define the closure of $R$ as 
    \[\cl(R) \coloneqq \bigcap_{G\in \mathfrak{G}} \rel_G.\]
\end{definition}

In particular, the last definition implies that $\cl(R)$ contains all pairs $(ab,xy)$ for which $\lca_G(ab)\preceq_G\lca_G(xy)$ 
must hold in any DAG $G$ that realizes $R$.

Now recall that in order to show that $\cl(R)$ 
is a ``valid'' closure for realizable relations, we must 
ensure that $\cl(R)$ satisfies the following three properties:
\begin{enumerate}[noitemsep]
  \item \emph{Extensivity:} $R \subseteq \cl(R)$.
  \item \emph{Monotonicity:} If $R_1 \subseteq R_2$ for two realizable relations on $\pairs(X)$, then $\cl(R_1) \subseteq \cl(R_2)$.
  \item \emph{Idempotency:} $\cl(\cl(R)) = \cl(R)$.
\end{enumerate}
A particular difficulty in proving these axioms directly lies in the fact that $\cl(R)$ is defined
only on \emph{realizable} relations. Hence, to show that $\cl(\cl(R))$ can be applied, we
must verify that $\cl(R)$ is realizable by some DAG. This, however, is not trivial. 
Clearly, $\cl(R)\subseteq \rel_G$ for all DAGs $G\in\mathfrak{G}$ that realize $R$. 
However, as Example~\ref{empl:subset-non-realizing} shows, this does not necessarily imply
that any of the DAGs $G\in\mathfrak{G}$ realize $\cl(R)$. In fact, as shown in 
Example~\ref{exmpl:not-subset-real}, there are subsets of realizable relations
that are not realizable. Thus, we cannot assume 
a priori that $\cl(R)$ is realizable by any DAG. 
To get around this issue, we could have defined $\cl(R)$ for all relations, realizable or not.
In this case, $\cl(R)=\emptyset$ for all non-realizable relations $R$. 
However, we would be faced with a similar problem. 
If $R$ is non-empty and realizable but $\cl(R)$ is not realizable, 
then $R\not\subseteq \cl(R)=\emptyset$; violating extensivity.

We now proceed with proving that $\cl(R)=R^+$ and thus, that $\cl(R)$ is realizable
for any realizable relation $R$.
We first establish some further results. The first 
(in particular Condition (2)) strengthens Lemma~\ref{lem:soundness}.

\begin{lemma}\label{lem:canonical-G-iff} 
	Let $R$ be a realizable relation on $\pairs(X)$ and $H\in \{\cG_R, \cN_R\}$ with $\cG_R$ being the canonical DAG and $\cN_R$ being the canonical network of $R$. 
    Then, the following statements hold:
    \begin{enumerate}
        \item For all $ab,xy\in \support_R^+$, $\lca_{H}(ab)$ and $\lca_{H}(xy)$ are well-defined, and the following two conditions hold
                \begin{enumerate}
                    \item $\lca_{H}(ab)\preceq_{H}\lca_{H}(xy) \iff ab\ R^+\ xy$.
                    \item $\lca_{H}(ab)\prec_{H}\lca_{H}(xy) \iff ab\ R^+\ xy \text{ and }(xy,ab)\notin R^+$.
                \end{enumerate}
        \item $H$ realizes $R^+$.                 
     \end{enumerate}
  \end{lemma}
  \begin{proof} 
  Let $R$ be a realizable relation on $\pairs(X)$ and $H\in \{\cG_R, \cN_R\}$ with $\cG_R$ being the canonical DAG and $\cN_R$ being the canonical network of $R$. 

    We start with proving Condition (1). Since $R$ is realizable, 
    \REV{ Theorem~\ref{thm:char} and Proposition~\ref{prop:srel-realized-by-N} imply that  $H$ realizes $R$.}
    By Lemma~\ref{lem:rel-subset}, $\lca_{H}(uv)$ is well-defined for all  $uv\in \support_R^+$. 
    We continue by proving the equivalence in Condition (1.a). 
    Let $ab,xy\in \support_R^+$.  The \emph{if}-direction is ensured
    by Lemma~\ref{lem:soundness}. For the \emph{only if}-direction, 
    assume that $\lca_{H}(ab)\preceq_{H}\lca_{H}(xy)$. 
    Recall that the vertices of $H$ consist of the leaves in $X$, the equivalence classes $[cd]$ with $cd\in\support_R$ and $c\neq d$ and, in case that $H=\cN_R$, the unique root $\rho$ that may or may not correspond to one of the
    equivalence classes (artificially added to
    \REVII{$\cG_R$}
    in the latter case). Since, by assumption, $ab,xy\in\support_R^+$, the equivalence classes $[ab]$ and $[xy]$ will by construction ``correspond'' to some vertex of $H$, but are possibly relabeled by a leaf in $X$. To circumvent this technicality, we
    let $v=[ab]$ if $a\neq b$ and, otherwise, $v=a$. Moreover, 
    let $w=[xy]$ if $x\neq y$ and, otherwise, $w=x$. 
    \REVII{Since $H\in \{\cG_R,\cN_R\}$, }
    Proposition~\ref{prop:CanonicalG-DAG-outsourced} respectively
    Proposition~\ref{prop:srel-realized-by-N} implies that 
    $v = \lca_{H}(ab)$	and $w = \lca_H(xy)$. 
    Since we have assumed that $\lca_{H}(ab)\preceq_{H}\lca_{H}(xy)$ and for \REVII{the two} possibilities of $H$ we have shown that $v = \lca_{H}(ab)$ and $\lca_{H}(xy) =w$, we conclude that $v\preceq_H w$.
	Clearly, if $H=\cG_R$, we have $v \preceq_{\cG_R}  w$. 
	If $H=\cN_R$, then the latter together with 
	the fact that \REVII{$\cG_R$} differs from $\cN_R$ only by possibly a newly added root 
	implies  $v \preceq_{\cG_R}  w$. In summary, $v\preceq_H w$ implies that $v\preceq_{\cG_R} w$, for any choice  $H \in  \{\cG_R, \cN_R\}$. 
	By construction of $\cG_R$, the latter ensures that $[ab]\le_{R^+} [xy]$
    and by definition of $\le_{R^+}$, $ab\ R^+\ xy$ must hold. 
    Therefore, Condition (1.a) holds. Condition (1.b) immediately follows from (1.a), 
    since $\lca_H(ab)\prec_H\lca_N(xy)$ if and only if $\lca_H(ab)\preceq_H\lca_H(xy)$ and $\lca_H(xy)\not\preceq_H\lca_H(ab)$. In summary, Condition (1) holds.

    We continue with proving Condition (2), which we have already proven for the case
    $H=\cG_R$ in Proposition~\ref{prop:CanonicalG-DAG-outsourced}. Since $R$ is realizable, we 
    can assume that Condition (1) holds. Let \REVII{$H=\cN_R$}.
    Since $R^+$ is transitive, it holds
    that $R^+ = \tc(R^+)$. Let $(ab,xy)\in R^+$. 
    If $(xy,ab)\notin \tc(R^+)=R^+$, 
    then  Condition (1.b) implies 
    that $\lca_H(ab)\prec_H \lca_H(xy)$, i.e., \axiom{I1} holds for $R^+$ and $H$. 
    If $(xy,ab)\in \tc(R^+)=R^+$,  then Condition (1.a) implies that 
    $\lca_H(ab)=\lca_H(xy)$, 
    i.e., \axiom{I2} holds for $H$ and $R^+$. Hence, $H$ realizes $R^+$.
  \end{proof}

  Lemma~\ref{lem:canonical-G-iff} shows that the realizability of $R$ implies 
  realizability of $R^+$. However, $R$ might be realized by some DAG $G$ that
  does not realize $R^+$, i.e., Conditions (1) and (2) in Lemma~\ref{lem:canonical-G-iff}  are not necessarily satisfied for arbitrary DAGs. To see this, we provide an example.
  
  \begin{example}[$G$ realizes $R$ but not $R^+$]\label{exmpl:GrealRnotR+}
  Consider the relation $R= \{(aa,ac), (bb,ac),(aa,ab)\}$ on \REV{$\pairs(X)$ with} $X=\{a,b,c\}$
  as shown in Example~\ref{exmpl-not-minimal}. As argued in Example~\ref{exmpl-not-minimal}, 
  the star-tree
  $T$ on $X$ realizes $R$ and we have $(ab,ac)\in R^+$ and $(ac,ab)\notin R^+$.
  Since $R^+$ is transitive, we have $\tc(R^+)=R^+$. 
  Hence, $(ab,ac)\in R^+$ and $(ac,ab)\notin R^+$ implies together with
  \axiom{I1} that any DAG $G$ realizing $R^+$ must satisfy $\lca_G(ab)\prec_G \lca_G(ac)$. 
  Since $\lca_T(ab) = \lca_T(ac)$ it follows that $T$ does not realize $R^+$. 
  Hence,  Condition (2) in Lemma~\ref{lem:canonical-G-iff}  does not hold for $H=T$. 
  In particular, $\lca_T(ac) \preceq_T \lca_T(ab)$ holds, but
  $(ac,ab) \notin R^+$, i.e., Condition (1) in Lemma~\ref{lem:canonical-G-iff}  
  does not hold. 
\end{example}

As stated in Section~\ref{sec:closure}, $R^+$ contains additional 
information about any DAG~$G$ that could potentially realize $R$. 
Example~\ref{exmpl:GrealRnotR+} may seem somewhat in contrast to this statement, 
since $G$ may realize $R$ but does not necessarily realize $R^+$. 
However, in light of Lemma~\ref{lem:soundness}, the statement is consistent: 
for any DAG~$G$ realizing $R$, if
$ab \ R^+ \ xy$ then $\lca_G(ab) \preceq_G \lca_G(xy)$. 
In particular, for Example~\ref{exmpl:GrealRnotR+},
we have $(ab,ac) \in R^+$ and 
$\lca_T(ab) \preceq_T \lca_T(ac)$ for the star-tree realizing
$R = \{(aa,ac), (bb,ac), (aa,ab)\}$. 
We can go even further: we now show
that for every realizable relation $R$, there exists a network~$G$ 
that realizes $R$ and satisfies $R^+ = \rel_G$, which strengthens 
the fact that $R^+\subseteq \rel_G$ (cf.\ Lemma~\ref{lem:soundness}).

\begin{lemma}\label{lem:G-st-relG=R+}
	For every realizable relation $R$, there exists a network $G$ that realizes $R$ and
	satisfies $\rel_G=R^+$. 
\end{lemma}
\begin{proof}
    Let $R$ be a realizable relation on $\pairs(X)$. By Lemma~\ref{lem:canonical-G-iff}, $R^+$ is realized by the canonical network $N\coloneqq \cN_R$. Recall from Theorem~\ref{thm:rules} that $\support_{R^+}=\support_R^+$.
    Let $\rho_N$ denote the unique root of $N$.
  In the following, we say that a DAG $G$ with $V(N)\subseteq V(G)$ is $\preceq_N$-preserving if the following property is satisfied: $v \preceq_N u$ if and only if $v \preceq_G u$ for all $u,v  \in V(N)$.
  Note that, by definition, $xx\in \support_{R}^+$ for all $x\in X$. Therefore, $ab \notin \support_{R^+}=\support_{R}^+$ implies $a\neq b$.

  Now, let $G$ be the directed graph obtained from $N$ by adding, 
  for all $ab \in \pairs(X)$ such that $ab \notin \support_{R^+}$, two new vertices
  $v_{ab}$ and $u_{ab}$, and the six arcs $v_{ab} \to a$, $v_{ab} \to b$, $\rho_N\to v_{ab}$, $u_{ab} \to a$, $u_{ab} \to b$ and $\rho_N\to u_{ab}$.
  As $N$ is a network on $X$, one easily observes that the leaf set of $G$ remains $X$. In $G$, both the vertex $v_{ab}$ and the vertex $u_{ab}$ have $\rho_N$ as its unique parent and precisely two children, and these two children are leaves in $X$. 
  It is, therefore, straightforward to verify that $G$ remains a network.
  By construction, $V(N) \subseteq V(G)$ and $G$ is $\preceq_N$-preserving.

  We now show that $\rel_G=R^+$.
  Suppose first that $ab,xy \in \pairs(X)$ are such that $ab\ R^+\ xy$. Since $N$ realizes $R^+$,
  the LCAs $\lca_N(ab)$ and $\lca_N(xy)$ must be well-defined, and satisfy $\lca_N(ab) \preceq_N
  \lca_N(xy)$. Moreover, $ab\ R^+\ xy$ implies that $ab,xy \in \support_{R^+}$. 
  Note that if $v\in V(G) \setminus V(N)$, then $v=v_{cd}$ or $v=u_{cd}$ for some $cd\notin \support_{R^+}$.
  In particular, the descendants of $v_{cd}$ are precisely $v_{cd}$, $c$ and $d$ and the descendants of $u_{cd}$ are precisely $u_{cd}$, $c$ and $d$. 
  Taking the latter arguments together, 
  there is no vertex $v \in V(G) \setminus V(N)$ such that $v$
  is an ancestor of both $a$ and $b$ in $G$, and there is no vertex $v \in V(G) \setminus V(N)$ such
  that $v$ is an ancestor of both $x$ and $y$ in $G$. 
  In other words, all common ancestors of $a$ and $b$ in $G$ as well as of $x$ and $y$ in $G$ are located in $V(N)\subseteq V(G)$.
  This, together with the fact that $G$ is $\preceq_N$-preserving, implies that $\lca_G(ab)=\lca_N(ab)$ and
  $\lca_G(xy)=\lca_N(xy)$ are well-defined and, moreover,  that $\lca_G(ab) \preceq_G \lca_G(xy)$.
  To summarize up to this point, we have shown that $ab\ R^+\ xy$ implies that $ab\rel_G xy$ i.e. that $R^+\subseteq \rel_G$.

   To show that $\rel_G\subseteq R^+$, suppose that $ab,xy \in \pairs(X)$ are such that $\lca_G(ab)$ and $\lca_G(xy)$ are well-defined, 
   and satisfy $\lca_G(ab) \preceq_G \lca_G(xy)$. We argue first that 
   $ab \in \support_{R^+}$. To this end, assume for contradiction, 
   that $ab \not\in \support_{R^+}$. In this case, the new vertices $v_{ab}$ and $u_{ab}$
   were added to $G$. It is easy to verify that,
   by construction, $v_{ab}$ and $u_{ab}$ are both contained in $\LCA_G(ab)$, which
   contradicts the assumption that $\lca_G(ab)$ is well-defined. Hence,  $ab \in \support_{R^+}$ must hold.
   Using similar arguments, we can conclude that $xy \in \support_{R^+}$.
   Since $\lca_G(ab)$ is, by assumption, well-defined and by the arguments in the preceding paragraph, we have 
   $\lca_G(ab)\in V(N)$. In particular, $\lca_G(ab)=\lca_N(ab)$ holds as $G$ is  $\preceq_N$-preserving. 
   Similarly,  $\lca_G(xy)=\lca_N(xy)$.
     Since $\lca_G(ab) \preceq_G
   \lca_G(xy)$ and since $G$ is $\preceq_N$-preserving it follows that $\lca_{N}(ab) \preceq_N
   \lca_{N}(xy)$. This together with Lemma~\ref{lem:canonical-G-iff} implies that $ab\ R^+\ xy$.
   We have thus shown that $\rel_G\subseteq R^+$, concluding the equality $\rel_G= R^+$. 
   
   Finally, we show that $G$ realizes $R$.
   To this end, let \REVII{$(ab,cd)\in R$}. 
   \REVII{Since $R\subseteq R^+=\rel_G$} it follows that $\lca_G(ab)$ and $\lca_G(cd)$
   are well-defined and satisfy $\lca_G(ab)\preceq_G \lca_G(cd)$. 

    There are now two cases to consider, either $(cd,ab)\notin\tc(R)$ or $(cd,ab)\in\tc(R)$.
    Suppose first that $(cd,ab)\notin\tc(R)$. Since $R$ is realizable, Theorem~\ref{thm:char} ensures that $R$
    satisfies \axiom{X2}. Hence, $(cd,ab)\notin R^+$. 
    This and $\rel_G= R^+$ implies that $\lca_G(cd)\not\preceq_G\lca_G(ab)$. 
    This together with  $\lca_G(ab)\preceq_G \lca_G(cd)$ implies that $\lca_G(ab)\prec_G \lca_G(cd)$. 
    Thus, Condition \axiom{I1} is satisfied. 
    Suppose now that  $(cd,ab)\in\tc(R)$. Since, by Observation~\ref{obs:tcR-subset-R+},
    $\tc(R)\subseteq R^+$ and, by similar arguments as before, 
    $\lca_G(cd)\preceq_G \lca_G(ab)$. This and $\lca_G(ab)\preceq_G \lca_G(cd)$ implies that
    $\lca_G(ab) = \lca_G(cd)$. Thus, Condition \axiom{I2} is satisfied. 
    In summary, $G$ realizes $R$.
\end{proof}

We now prove that $\cl(R)=R^+$ and derive some further properties of $\cl(R)$. 

\begin{theorem}\label{thm:closure-classic}
   	   	Let  $R$ be a realizable relation on $\pairs(X)$. Then, $\cl(R)=R^+$ and $\cl(R)$ is realizable.
    Furthermore, $\cl$ is a closure operator, i.e., it 
    satisfies the classical closure axioms \emph{Extensivity}, \emph{Monotonicity}, and \emph{Idempotency}. 
	In addition,  $\cl(R)$ can be computed in polynomial time in $|X|$ \REV{as specified in Algorithm~\ref{alg:Rplus} }.  
\end{theorem}
\begin{proof}
	Let  $R$ be as in the statement of the theorem. Furthermore, 
    let  $\mathfrak{G}$ denote the set of all DAGs $G$ on $X$ that realize $R$
    and  $\mathfrak{R}$ the set of all relations $S$ on $\pairs(X)$ 
    that are $\support_R^+$-reflexive, transitive, cross-consistent, and satisfy $R \subseteq S$. 

    Let $G\in \mathfrak{G}$. By Proposition~\ref{prop:G-rel-relplus}, $\rel_G^+ = \rel_G$. 
	Hence, $\rel_G$ is  transitive, $\support_{\rel_G}^+$-reflexive and cross-consistent. 
	  Lemma~\ref{lem:rel-subset} implies that
    $R \subseteq \rel_G$ which in particular ensures that $\support_R\subseteq\support_{\rel_G}$. As $G$ is a DAG on $X$, $\rel_G$ is a relation on $\pairs(X)$ and, thus, $\support_R^+\subseteq\support_{\rel_G}^+$. The latter two facts together with $\rel_G$ being $\support_{\rel_G}^+$-reflexive ensure that $\rel_G$ is $\support_R^+$-reflexive. Together with transitivity and cross-consistency of $\rel_G$, we conclude that $\rel_G\in \mathfrak{R}$.
    Hence, for all $G\in \mathfrak{G}$ we have $\rel_G\in \mathfrak{R}$
    which implies that $R^+  =\cap_{R'\in\mathfrak{R}} R'
    \subseteq \cap_{G\in \mathfrak{G}} \rel_G = \cl(R)$. 
    Moreover, since $R$ is realizable, 
    Lemma~\ref{lem:G-st-relG=R+} ensures the existence of a DAG $G\in\mathfrak{G}$ for which $\rel_G=R^+$. By definition, $\cl(R)\subseteq\rel_G=R^+$. 
    Thus, $\cl(R)=R^+$.

Since $R$ is realizable, it satisfies \axiom{X1} (cf.\ Theorem~\ref{thm:char}). 
    By Proposition~\ref{prop:CanonicalG-DAG-outsourced}, the canonical DAG $\cG_R$ realizes
    $R^+=\cl(R)$. Hence, $\cl(R)$ is realizable.
    
    Furthermore, since $\cl(R)=R^+$ and 
    since $R^+$ satisfies the closure axioms  (cf.\ Proposition~\ref{lem:closure-axioms})
    if follows that $\cl$ satisfies the closure axioms. Moreover, $\cl(R)=R^+$
    together with Theorem~\ref{thm:rules} \REV{and \ref{thm:R+-polytime}} imply that 
    $\cl(R)$ can be computed in polynomial time in $|X|$ by \REV{using Algorithm~\ref{alg:Rplus}}.
\end{proof}

Observation~\ref{obs:G_R=G_R+} together with Theorem~\ref{thm:closure-classic} immediately implies:
\begin{corollary}
    For all realizable relations $R$ we have $\cG_R=\cG_{\cl(R)}$.
\end{corollary}

Intriguingly, using the classical closure, we now show that 
for all DAGs $G$, the information contained in $\rel_G$ is entirely determined
by that contained in $\srel_G$ and its $+$-closure $\srel_G^+$.

\begin{proposition}\label{prop:cor:G-rel-relplus}
For all DAGs $G$ it holds that $\rel_G=\rel_G^+ = \cl(\rel_G)= \srel_G^+ = \cl(\srel_G)$ and $\srel_G\neq \srel_G^+$. In particular, $\support_{\rel_G}^+ = \support_{\srel_G}^+$.
Moreover, for all DAGs $G$ and $H$ \mh{on the same leaf set $X$}, it holds that  
    \[\srel_G=\srel_H \iff \rel_G=\rel_H.\] 
\end{proposition}
\begin{proof}
    Let $G$ be a DAG. 
    \mh{By Lemma~\ref{lem:G-realizes-relG}, $\rel_G$ is realizable 
    By Proposition~\ref{prop:G-rel-relplus} and Theorem~\ref{thm:closure-classic}, 
    we have $\rel_G=\rel_G^+ = \cl(\rel_G)$. Moreover, 
    by Lemma~\ref{lem:G-realizes-relG}, 
    $\srel_G$ is strictly-realizable. This and 
    Lemma~\ref{lem:str-real=>real} implies that 
    $\srel_G$ is realizable. Theorem~\ref{thm:closure-classic}
    implies that $\srel_G^+ = \cl(\srel_G)$.
    }
    
    Hence, it suffices to show that $\rel_G = \srel_G^+$. By definition, 
    $\srel_G\subseteq \rel_G$. By Proposition~\ref{lem:closure-axioms}, 
    $\srel_G^+\subseteq \rel_G^+$. As $\rel_G^+=\rel_G$, we have
    $\srel_G^+\subseteq \rel_G$. 
    To show that $\rel_G \subseteq \srel_G^+$, let $(ab,xy)\in \rel_G$
    and thus, $\lca_G(ab)$ and $\lca_G(xy)$ are well-defined and
    satisfy  $\lca_G(ab)\preceq_G\lca_G(xy)$. 
    If $\lca_G(ab)\prec_G\lca_G(xy)$, then $(ab,xy)\in \srel_G\subseteq \srel_G^+$
    and we are done. Suppose that  $\lca_G(ab) = \lca_G(xy)$. 
    If $a=b$, then $a = \lca_G(ab) = \lca_G(xy)$ implies that $x=y=a$. 
    Since $\srel_G^+$ is $\support_{\srel_G}^+$-reflexive, 
    we have $(aa,aa) = (ab,xy)\in \srel_G^+$. 
    Suppose that $a\neq b$. Note that, in this case, 
    $\lca_G(aa)\prec_G\lca_G(ab)$ and, therefore, $ab\in \support_{\srel_G}
    \subseteq \support_{\srel_G^+}$. 
    Since $\lca_G(ab) = \lca_G(xy)$, we have 
    $aa\ \srel_G\ xy$ and $bb\ \srel_G\ xy$. 
    This together with $ab\in\support_{\srel_G^+}$ and cross-consistency of $\srel_G^+$ 
    implies that $(ab,xy)\in \srel_G^+$.
    Hence, $\rel_G \subseteq \srel_G^+$ and, therefore, $\rel_G =\srel_G^+$
    holds. In summary, $\rel_G=\rel_G^+ = \cl(\rel_G)= \srel_G^+ = \cl(\srel_G)$. 
    This together with Theorem~\ref{thm:rules} implies that $
    \support_{\srel_G}^+=\support_{\srel_G^+} = \support_{\rel_G^+}=\support_{\rel_G}^+$. 

    We show now that $\srel_G\neq \srel_G^+$. 
    Observe first that $\srel_G\neq \rel_G$, 
    since for any leaf $x\in X$, we have $(xx,xx)\in \rel_G$ but $(xx,xx)\notin \srel_G$.
    As shown above, $\srel_G^+ = \rel_G$. 
    Taking the latter two arguments together, we obtain $\srel_G\neq \srel_G^+$.

    Finally, let $G$ and $H$ be two DAGs \mh{on the same leaf set $X$}. We show that 
    $\srel_G=\srel_H \iff \rel_G=\rel_H$.
    \mh{If $\srel_G=\srel_H$, then $\support_{\srel_G}^+=\support_{\srel_H}^+$ since $G$ and $H$ have the same set of leaves. Therefore, 
    $\srel_G=\srel_H$ implies $\srel_G^+=\srel^+_H$. As previously argued, it holds that $\rel_G=\srel^+_G = \srel^+_H = \rel_H$.}
    Suppose now that $\rel_G=\rel_H$. Assume, for contradiction that $\srel_G\neq \srel_H$. 
    Thus, we can, without loss of generality, assume that there is some $(ab,xy)\in \srel_G$ such that 
    $(ab,xy)\notin \srel_H$. Since $(ab,xy)\in \srel_G$, we have $\lca_G(ab)\prec_G \lca_G(xy)$.
    Since $\srel_G\subseteq \rel_G$ and  $\rel_G=\rel_H$, 
    if follows that $(ab,xy)\in \rel_H$ and, hence, that 
    $\lca_H(ab)\preceq_H\lca_H(xy)$. 
    This together with $(ab,xy)\notin \srel_H$ implies 
    $\lca_H(ab)=\lca_H(xy)$. Hence, 
    $(xy,ab)\in \rel_H$. Since $\rel_H=\rel_G$ it follows that 
    $(xy,ab)\in \rel_G$ and, therefore, 
    $\lca_G(xy)\preceq_G\lca_G(ab)$; a contradiction to $\lca_G(ab)\prec_G \lca_G(xy)$. 
    Hence, $\srel_G = \srel_H$. 
\end{proof}

Interestingly, using this last result, we 
can also show that it is sufficient to know $\rel_H$ and $\support_R^+$ to determine $\cl(R)$ 
for any realizable relation $R$, where $H$ is its canonical DAG or network.

\begin{theorem}
Let $R$ be a realizable relation and  $H\in \{\cG_R, \cN_R\}$ with $\cG_R$ being the canonical DAG
    and $\cN_R$ being the canonical network of $R$. Then it holds that 
\[\cl(R) = \rel_{H} \cap (\support_R^+ \times \support_R^+).\]
\end{theorem}
\begin{proof}
Let $R$ be a realizable relation on $\pairs(X)$ and $H\in \{\cG_R,
\cN_R\}$ be one of the DAGs
    on $X$. For both choices of $H$, Theorem~\ref{thm:char} respectively Proposition~\ref{prop:srel-realized-by-N} implies that $H$ realizes $R$. 
   This and Lemma~\ref{lem:soundness} implies that $R^+\subseteq\rel_H$.
   By Theorem~\ref{thm:closure-classic}, $\cl(R) = R^+$ and, therefore, 
    $\cl(R)\subseteq\rel_H$. Moreover, $\cl(R) = R^+$ together with Theorem~\ref{thm:rules} implies
    that $\support_R^+ =\support_{R^+} = \support_{\cl(R)}$. Taking the latter two arguments
    together with the fact that $\cl(R) = \cl(R) \cap (\support_{\cl(R)} \times \support_{\cl(R)})$,
    we obtain \[\cl(R) = \cl(R) \cap (\support_R^+ \times \support_R^+) \subseteq
    \rel_{H}\cap (\support_R^+ \times \support_R^+) .\] 

    Now, let $(ab,xy)\in \rel_{H}\cap(\support_R^+ \times \support_R^+)$. By definition of $\rel_H$ we have  $\lca_{H}(ab)\preceq_{H} \lca_{H}(xy)$, which taken together with $ab,xy\in \support_R^+$ and Lemma~\ref{lem:canonical-G-iff}(1.a) ensures that $(ab,xy)\in R^+$. We therefore have that $ \rel_{H}\cap(\support_R^+ \times \support_R^+)\subseteq R^+=\cl(R)$, where the last equality is ensured by Theorem~\ref{thm:closure-classic}. This concludes that $\cl(R) =
    \rel_{H}\cap (\support_R^+ \times \support_R^+)$.
\end{proof}

We conclude this section by considering the relations that remain invariant under the 
classical closure operator. In particular, we want to characterize which 
realizable relations are closed, i.e., those $R$ that satisfy  $R=\cl(R)$.  
By Theorem~\ref{thm:closure-classic} one 
characterization for such relations is that $R=R^+$ holds. 
In the following theorem, we characterize closed relations in 
terms of a ``stronger'' form of realizability. 

\begin{theorem}\label{thm:cl(R)-canonDAG}
    Let $R$ be a realizable relation. Then $R=\cl(R)$ if and only if there exists a DAG $H$ such that for all $ab,cd\in\support_R^+$, $\lca_H(ab)$ and $\lca_H(cd)$ are well-defined and satisfy
    \begin{equation}\label{eq:iff-R-lca}
        ab\ R\ cd \iff \lca_H(ab)\preceq_H\lca_H(cd).
    \end{equation}
    In particular, $R=\cl(R)$ if and only if the Equivalence~\eqref{eq:iff-R-lca}
    holds for $R$ and  $H\in \{\cG_R, \cN_R\}$ with $\cG_R$ being the canonical DAG
    and $\cN_R$ being the canonical network of $R$.
\end{theorem}
\begin{proof}
    Let $R$ be a realizable relation on $\pairs(X)$. Recall from Theorem~\ref{thm:closure-classic} that $\cl(R) = R^{+}$, 
    a fact that will be used repeatedly throughout this proof.
    The \emph{only if}-direction is an immediate consequence of Lemma~\ref{lem:canonical-G-iff} 
    since, if $R=\cl(R)=R^+$, then $ab\ R^+\ xy$ holds if and only if $ab\ R\ xy$, for all $ab,xy\in\support_R^+$.

   For the \emph{if}-direction, assume that there is a DAG $H$ such that, for all $ab,cd\in\support_R^+$, $\lca_H(ab)$ and $\lca_H(cd)$ are well-defined and the Equivalence in \eqref{eq:iff-R-lca} is satisfied.
     We proceed by showing that $R$ is $\support_R^+$-reflexive, transitive and cross-consistent. Since $\lca_H(ab)\preceq_H\lca_H(ab)$ for all $ab\in\pairs(X)$ for which $\lca_H(ab)$ is well-defined, $\lca_H(ab)\preceq_H\lca_H(ab)$ in particular holds for all $ab\in\support_R^+$. By Eq.~\eqref{eq:iff-R-lca}, it follows that $R$ is $\support_R^+$-reflexive. 
    If $ab\ R\ xy$ and $xy\ R\ cd$ holds, Eq.~\eqref{eq:iff-R-lca} ensures that $\lca_H(ab)\preceq_H\lca_H(xy)$ and $\lca_H(xy)\preceq_H\lca_H(cd)$ which, by transitivity of $\preceq_H$, shows that $\lca_H(ab)\preceq_H\lca_H(cd)$. Hence, Eq.~\eqref{eq:iff-R-lca} implies that $ab\ R\ cd$ i.e. $R$ is transitive. Finally, suppose $ab\ R\ xy$, $cd\ R\ xy$ and $ac\in\support_R$. Thus, the least common ancestors $\lca_H(ac)$, $\lca_H(ab)$, $\lca_H(cd)$ and $\lca_H(xy)$ are well-defined and satisfy, by Eq.~\eqref{eq:iff-R-lca}, $\lca_H(ab)\preceq_H\lca_H(xy)$ and respectively $\lca_H(cd)\preceq_H\lca_H(xy)$. Lemma~\ref{lem:simpleLCAinfer} implies  $\lca_H(ac)\preceq_H\lca_H(xy)$ which together with Eq.~\eqref{eq:iff-R-lca} and the fact that $ac\in\support_R\subseteq \support_{R}^+$
    implies that $ac\ R\ xy$. In other words, $R$ is cross-consistent. Since we have shown that $R$ is $\support_R^+$-reflexive, transitive and cross-consistent, Theorem~\ref{thm:rules} ensures that $R=R^+=\cl(R)$. 

    We prove now the last statement in the theorem.
    Suppose that $R=\cl(R)=R^+$. 
        Since $R$ is realizable, Theorem~\ref{thm:char} respectively
        Proposition~\ref{prop:srel-realized-by-N} implies that 
        $H\in \{\cG_R, \cN_R\}$ realizes $R$.
        Moreover, Lemma~\ref{lem:canonical-G-iff}(1) implies that, for all
        $ab,cd\in\support_R^+$, it holds that 
        $ab\ R^+\ cd \iff \lca_H(ab)\preceq_H\lca_H(cd)$. Since $R=R^+$, 
        the Equivalence~\eqref{eq:iff-R-lca} holds for $R$ and  $H\in \{\cG_R, \cN_R\}$. 
        Suppose now that  $ab\ R\ cd \iff \lca_H(ab)\preceq_H\lca_H(cd)$
        holds for $R$ and $H\in \{\cG_R, \cN_R\}$. 
        By Lemma~\ref{lem:canonical-G-iff}(1), we have 
        $ab\ R^+\ cd \iff \lca_H(ab)\preceq_H\lca_H(cd)$.
        It now immediately follows that $R=R^+=\cl(R)$.  
\end{proof}


\section{Involving Incomparability Constraints} 
\label{sec:incomp}

Given a collection of taxa $X$ (such as species or genes) and three taxa $a, b, c \in X$, 
in a phylogenetic tree $T$ on $X$, we say that $a$ is evolutionary more closely related to $b$ 
than to $c$ if $\lca_T(ab) \prec_T \lca_T(ac)$. 
This holds because, for any vertex $v$ in $T$, the vertex $\lca_T(av)$ always lies on the 
(unique) path between the root of $T$ and $v$. 
Hence, if $\lca_T(ab) \prec_T \lca_T(ac)$, the ``branching point'' corresponding to $\lca_T(ab)$ 
must have occurred more recently (i.e., closer to the present) than the branching point 
corresponding to $\lca_T(ac)$. 
This notion of evolutionary relatedness has proven useful in the characterization and inference 
of so-called best matches~\cite{Stadler2020,Geiss2019,Geiss2020}.

\begin{figure}
    \centering
    \includegraphics[width=0.8\textwidth]{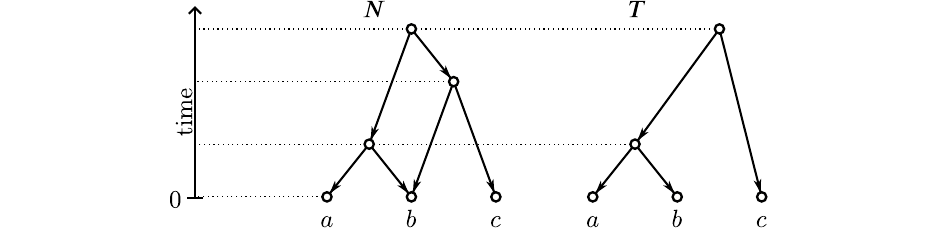}
    \caption{A phylogenetic tree $T$ and network $N$ on the set $X=\{a,b,c\}$ with time from the present indicated. In both DAGs, $\lca(ab)$, $\lca(ac)$ and  $\lca(bc)$ are well-defined, and
             $a$ is evolutionary closer to $b$ than to $c$. The latter holds because
             the time elapsed since $\lca(ab)$ existed is
             less than that since $\lca(ac)$ existed. However,
             $\lca_N(ab)$ and $\lca_N(ac)$ are $\preceq_N$-incomparable.}
    \label{fig:evol}
\end{figure}

This definition of evolutionary relatedness, however, does not directly translate to phylogenetic networks. 
In particular, it may occur that the vertices $\lca_N(ab)$ and $\lca_N(ac)$ are both well-defined but 
$\preceq_N$-incomparable, even though $a$ is still evolutionary closer to $b$ than to $c$; 
see Figure~\ref{fig:evol} for an illustrative example.
Such situations naturally arise in reticulate evolutionary processes such as hybridization or horizontal gene transfer, 
where different parts of the genome may follow distinct ancestral histories. 
It is therefore natural to introduce \emph{incomparability constraints}, capturing pairs of LCA vertices that are not ordered by $\preceq_N$. 
These constraints provide both a combinatorial means to represent structural incomparability in the ancestral relation 
and a biologically meaningful way to describe lineages that coexist without a clear hierarchical relationship.

To capture these ideas more formally we make the following definition.

\begin{definition}\label{def:pair-realize}
  Let $R$ and $S$ be two relations on $\pairs(X)$. 
  We say that the ordered pair $(R,S)$ is \emph{realized} by a DAG $G$ 
  on $X$ if the following two conditions hold.
  \begin{enumerate}
  	\item $G$ realizes $R$.
   \item For all $(ab, xy)\in S$ both $\lca_G(ab)$ and $\lca_G(xy)$ are well-defined and                 $\preceq_G$-incomparable. 
 \end{enumerate}
  In this case, we also say that $(R,S) $ is \emph{realizable}.
\end{definition}

In the rest of this section we shall derive a characterization for when
an ordered pair of relations is realizable (see Theorem~\ref{thm:incomp}).
To this end, we first consider when the union of two realizable relations is realizable.
Note that, as the following example shows, the union may not be realizable.

\begin{example}[Union of  realizable relations that is not realizable]\label{exmpl:union-not-realized}
Consider the two trees $T_3$ and $T_4$ as shown in Figure~\ref{fig:realize-subset-ex} and the
 two relations $R_1 =\{(xx,ab), (yy,ab)\}$ and 
 $R_2 =\{ (ab,xy)\}$ on $\pairs(X)$ with $X=\{a,b,x,y\}$.
  One can easily check that $R_1=\tc(R_1)$ and  $R_2=\tc(R_2)$
 and that $R_1$ is realized by the star-tree $T_4$ while $R_2$ is realized by the tree $T_3$. 
 However, $R = R_1\cup R_2$ is exactly the relation as  in Example~\ref{exmpl:not-subset-real}, where we have argued that $R$ is not realizable. Note that here $\tc(R_1)\cup \tc(R_2)\neq \tc(R_1\cup R_2)$
 since $(xx,xy)\in \tc(R)$ but $(xx,xy)\notin \tc(R_1)\cup\tc(R_2)$.
\end{example}

Even though we have seen that the union of two realizable relations $R$ and $R'$ may not be realizable,
we now show that this is the case when $R$ and $R'$ can be realized by the same DAG.

\begin{lemma}\label{lm:union}
Suppose that $R$ and $R'$ are relations on $\pairs(X)$ such that there exists a DAG $G$ on $X$ realizing both $R$ and $R'$. Then, $G$ realizes $R \cup R'$.
\end{lemma}
\begin{proof}
Let $R$ and $R'$ be relations on $\pairs(X)$ that are both realized by the DAG $G$ on $X$.
\REV{Note that $\lca_G(ab)$ is well-defined for every $ab\in\support_{\tilde R}^+$, since $\support_{\tilde R}^+=\support_R^+\cup\support_{R'}^+$ by definition.}
We show that $\tilde R \coloneqq R \cup R'$ satisfies \axiom{I1} and \axiom{I2}.

\REVII{Suppose first that $(ab,cd) \in \tilde R$ and $(cd,ab) \notin \tc(\tilde R)$. Without loss of generality, assume $(ab,cd)\in R$. Since $R,R'\subseteq \tilde R$, 
monotonicity of the transitive closure ensures that
$\tc(R)\cup\tc(R')\subseteq\tc(\tilde R)$. Thus $(cd,ab)\notin\tc(\tilde R)$ implies $(cd,ab)\notin\tc(R)$. The latter together with $(ab,cd)\in R$ and the fact that $G$ realizes $R$ implies, via \axiom{I1}, that $\lca_G(ab)\prec_G\lca_G(cd)$. In other words, $G$ and $\tilde R$ satisfy \axiom{I1}.} 

Suppose that  \REVII{$(ab,cd) \in \tilde R$} and $(cd,ab) \in \tc(\tilde R)$ holds. 
\REVII{Since $G$ realizes both $R$ and $R'$, $(ab,cd)\in\tilde R$ and Lemma~\ref{lem:rel-subset} implies that $\lca_G(ab)\preceq_G\lca_G(cd)$. Moreover, since $(cd,ab) \in \tc(\tilde R)$, there is a $(cd,ab)$-chain in $\tc(\tilde R)$ and thus 
some $p_0,\dots,p_k\in \support_{\tilde R}$, $k\geq 1$ such that $cd=p_0\ \tilde R\ p_1\ \tilde R\  \dots\ \tilde R\  p_{k-1}\ \tilde R\ p_k=ab$. In particular, $p_j\ R^j\ p_{j+1}$ for all $j \in \{0, \ldots, k-1\}$ 
and some $R^j\in \{R,R'\}$. This together with Lemma~\ref{lem:rel-subset} and the fact that
$G$ realizes both $R$ and $R'$ implies that $\lca_G(p_j)\preceq_G\lca_G(p_{j+1})$
for all $j \in \{0, \ldots, k-1\}$. By transitivity of the $\preceq_G$-relation, 
it follows that $\lca_G(cd) \preceq_G \lca_G(ab)$.}
Therefore,  $\lca_G(ab)=\lca_G(cd)$ and \axiom{I2} holds. In summary, $\tilde R$ is realized by $G$.
\end{proof}

Now, for two relations $R$ and $S$ on $\pairs(X)$ we let 
	\[R_S \coloneqq R \cup \{(aa,ab), (bb,ab) \mid ab\in \support_{S}\setminus \support^+_{R} \text{ and } a\neq b\}.\]

\begin{lemma}\label{lem:widetilde} 
    Let $R$ and $S$ be two relations on $\pairs(X)$.
	  The pair $(R,S)$ is realized by the DAG $G$ if and only if
	  the pair $(R_S,S)$ is realized  by the DAG $G$.
\end{lemma}
\begin{proof} 
    Let $R$ and $S$ be two relations on $\pairs(X)$.
    Since $R\subseteq R_S$ it follows that $\tc(R)\subseteq \tc(R_S)$
    and $\support_R^+\subseteq \support_{R_S}^+$.
    
    For the \emph{if}-direction suppose that $G$ realizes $(R_S,S)$. We need to 
    show that $(R,S)$ is also realized by $G$, that is, both \axiom{I1} and  \axiom{I2} hold.
    We first consider \axiom{I2}. To this end, suppose that \REVII{$(ab,cd)\in R$}
    and $(cd,ab)\in \tc(R)$. Since \REVII{$R\subseteq R_S$ and} $\tc(R)\subseteq \tc(R_S)$ it follows that
    \REVII{$(ab,cd)\in R_S$} and $(cd,ab)\in \tc(R_S)$. Since $G$ realizes $R_S$,  \axiom{I2}
    implies $\lca_G(cd)=\lca_G(ab)$. Hence, \axiom{I2}
    holds also for $R$ with respect to $G$.
    We now consider  \axiom{I1}.
    Assume that \REVII{$(ab,cd)\in R$} and $(cd,ab)\notin \tc(R)$.
    Therefore, $ab \neq cd$. 
    Assume for contradiction 
   that $(cd,ab)\in \tc(R_S)$. Then, there exists a sequence $p_1, \ldots,p_k$, $k \geq 2$ 
   of pairwise distinct elements of $\pairs(X)$ such that $p_1=cd$, $p_k=ab$, and $(p_j,p_{j+1}) \in R_S$ for all $j \in \{1, \ldots, k-1\}$. Since $(cd,ab)\notin \tc(R)$,  there must exist one such $j$ such that $(p_j,p_{j+1}) \in R_S \setminus R$. By definition of the set $R_S \setminus R$, this implies $p_j=xx$ and $p_{j+1}=xy$ for some distinct $x,y \in X$. 
   Since \REVII{$R_S$} is realizable, Theorem~\ref{thm:char} ensures that \REVII{$R_S$} satisfies \axiom{X1}, so there is no element $q \neq \REVII{xx}$ in $\pairs(X)$ such that \REVII{$(q,xx)\in R_S$. Note that if $j>1$, then $p_j=xx$ together with $(p_{j-1},p_j)=(p_{j-1},xx)\in R_S$ would yield a contradiction to the latter argument. 
   In other words, $p_j=xx$ force $j=1$ and,} therefore, $p_j=p_1=cd$ must hold, which implies $c=d$. 
   \REVII{In particular, we have $(ab,cc)=(ab,cd)\in R\subseteq R_S$.} However, we have assumed that $ab \neq cc$, 
   and so we obtain a contradiction to \REVII{$R_S$ satisfying} \axiom{X1}. \REVII{Consequently}, $(cd,ab)\notin \tc(R_S)$ must hold and we can conclude,  by similar arguments as for \axiom{I2}, that \axiom{I1} holds for $R$ with respect to $G$. 
   
   For the \emph{only if}-direction assume that $(R,S)$ is realized by the DAG $G$ on $X$.
   We first show that $R_S \setminus R$ is also realized by $G$. \REVII{Let $(p,q) \in R_S \setminus R$.} By definition of $R_S \setminus R$, there exists $a,b \in X$ distinct such that $p=aa$ and $q=ab$. It also follows directly from the definition that $(ab,aa) \notin \tc(R_S \setminus R)$. Moreover, $ab \in \support_S$, so since $G$ realizes $(R,S)$, $\lca_G(ab)$ is well defined. Clearly, $a=\lca_G(aa) \prec_G \lca_G(ab)$ always holds. \REVII{Hence \axiom{I1} holds for $R_S\setminus R$ with respect to $G$, and \axiom{I2} is vacuously satisfied.} This concludes the proof that $G$ realizes $R_S \setminus R$. 
   Since, in addition, $G$ realizes $R$ by assumption and since
   $R_S=R \cup (R_S \setminus R)$, 
   Lemma~\ref{lm:union} implies that $G$ realizes $R_S$. Hence, $G$ realizes $(R_S,S)$.
\end{proof}

We now prove the main result of this section.

\begin{theorem}\label{thm:incomp}
 Let $R$ and $S$ be two relations on $\pairs(X)$.
 The pair $(R,S)$ is realizable if and only if
 \begin{enumerate}[label=(\alph*)]
    \item $R_S$ is realizable, and 
    \item $S\subseteq I\coloneqq\{(p,q)\in\pairs(X)\times\pairs(X)\mid 
           (p,q), (q,p)\notin  R_S^+\}$.
  \end{enumerate}
  Moreover, 
  if $(R,S)$ is realizable, then $(R,S)$ is realized by $H\in \{\cG_{R_S}, \cN_{R_S}\}$ 
  with $\cG_{R_S}$ being the canonical DAG and $\cN_{R_S}$ being the canonical network of $R_S$. 
  In particular, verifying whether $(R,S)$ is realizable and, if so, constructing a DAG on $X$ that
  realizes $(R,S)$ can be done in polynomial time in $|X|$.
\end{theorem}
\begin{proof}
	Suppose the pair $(R,S)$ is realized by the DAG $G$ on $X$. 
	By Lemma~\ref{lem:widetilde}, $(R_S,S)$  is realized by $G$. 
    In particular, $R_S$ is realizable, i.e., Condition (a) holds.
	Now assume, for contradiction, that (b) does not hold, i.e. $S\not \subseteq I$. Hence, there is some $(ab, xy)\in
	S$ such that $ab\ R_S^+\ xy$ or $ xy\ R_S^+\ ab$. 
	Since $G$ realizes $(R_S,S)$ and $(ab, xy)\in S$, we must have  that
    $\lca_G(ab) $ and $\lca_G(xy)$ are $\preceq_G$-incomparable. 
    However, Lemma~\ref{lem:soundness} together with $ab\ R_S^+\ xy$ or $ xy\ R_S^+\ ab$
    implies that
	$\lca_G(ab)\ \preceq_G\ \lca_G(xy)$ or $\lca_G(xy)\ \preceq_G\ \lca_G(ab)$;
	a contradiction. Hence, Condition (b) holds. 

Conversely, suppose that $R_S$ is realizable and  $S\subseteq I$. 
By Theorem~\ref{thm:char} respectively Proposition~\ref{prop:srel-realized-by-N}, $R_S$ is realized by
 $H\in \{\cG_{R_S}, \cN_{R_S}\}$. 
 Now, take
    $ab,xy\in\pairs(X)$ such that $ab\ S\ xy$. By definition of $R_S$, we have
    $ab\in\support_{R_S}\subseteq \support_{R_S}^+$ and therefore
    $\lca_H(ab)$ is well-defined in $H$ (c.f.\ Lemma~\ref{lem:rel-subset}).  
    Analogously, $\lca_H(xy)$ is well-defined in $H$. Suppose, for contradiction, that $\lca_H(ab)$ and $\lca_H(xy)$
    are $\preceq_H$-comparable. We may, without loss of generality,  assume that
    $\lca_H(ab)\preceq_H\lca_H(xy)$. In this case, Lemma~\ref{lem:canonical-G-iff} implies that
    $ab\ R_S^+\ xy$ holds; a contradiction to $(ab,xy)\in S\subseteq I$. 
    Hence, $\lca_H(ab)$ and $\lca_H(xy)$ are $\preceq_N$-incomparable
    for all $(ab,xy)\in S$. Therefore, $H$ realizes the pair $(R_S,S)$.
    By Lemma~\ref{lem:widetilde}, $(R,S)$ is realized by $H$. 
    
    Finally, it is straightforward to see that checking if $S \subseteq I$ holds can be done in polynomial time in $|X|$. 
    By Theorem~\ref{thm:alg:char}, testing whether $R_S$ is realizable and, if so, constructing 
    $H$ can also be accomplished in polynomial time in $|X|$. 
    Together with the fact that $(R,S)$ is realizable if and only if Conditions~(a) and~(b) are satisfied 
    (in which case $(R,S)$ is realized by $H$ as argued above), 
     this implies that one can decide in polynomial time in $|X|$ whether $(R,S)$ is realizable and, if so, 
     construct a DAG on $X$ that realizes $(R,S)$. 
\end{proof}


\section{Discussion and Outlook} \label{sec:outlook}

In this paper, we have characterized when a set of LCA constraints is
realizable by a DAG or a phylogenetic network. Moreover, by relating realizability to
generalizations of classical 
closure operators for phylogenetic trees, we were able to develop a polynomial-time algorithm for 
deciding whether or not a set of LCA constraints is realizable, and
for constructing a DAG or network realizing this set if this is the case.
The algorithms introduced in this work have been implemented in Python and are freely available at \cite{github-AL}.
There remain a number of interesting future directions to explore.

\begin{figure}
    \centering
    \includegraphics[width=0.8\textwidth]{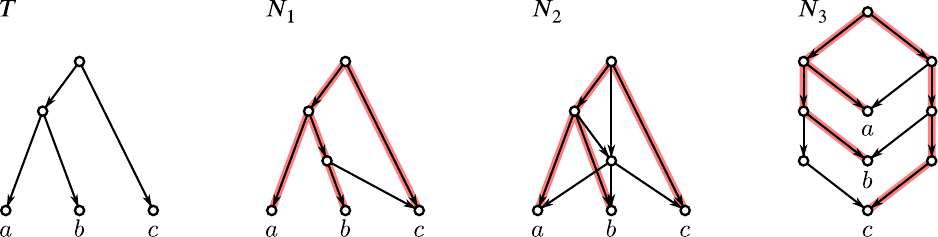}
    \caption{Four networks $T,N_1,N_2$ and $N_3$ on the set $X=\{a,b,c\}$. 
    The phylogenetic tree $T$ displays the triplet $ab|c$ according to \axiom{T1}
    and \axiom{T2}.
    All three networks $N_1$, $N_2$ and $N_3$
    display the triplet $ab|c$ according to \axiom{T2}
    (highlighted by the shaded red arcs)  but not according to \axiom{T1}. 
    In $N_1$ we have $\lca_{N_1}(bc)\prec_{N_1} \lca_{N_1}(ab)=\lca_{N_1}(ac)$. 
    Hence, instead of $ab|c$,  $N_1$ displays the triplet $bc|a$ according to \axiom{T1}. 
    In $N_2$, we have $\lca_{N_2}(bc)= \lca_{N_2}(ab)=\lca_{N_2}(ac)$. Hence, none of the triplets
    $ab|c$, $bc|a$ and $ac|b$ is displayed in $N_2$ according to \axiom{T1}. 
    The network $N_3$ displays all of the triplets $ab|c$, $bc|a$ and $ac|b$ according to \axiom{T2}. 
    However, none of the LCAs $\lca_{N_3}(ab)$, $\lca_{N_3}(bc)$ and $\lca_{N_3}(ac)$
    is well-defined in $N_3$, i.e., $N_3$ does not display any of triplets
    $ab|c$, $bc|a$ and $ac|b$ according to \axiom{T1}. }
    \label{fig:tiplesVSlca}
\end{figure}

\paragraph{Triplets and Networks.}
Approaches to reconstructing phylogenetic trees that employ \textsc{Build} usually rely on 
{\em triplets}, i.e., binary phylogenetic trees on three leaves (see the tree 
$T$ in Figure~\ref{fig:tiplesVSlca} for an illustrative example).
Similarly, the literature on phylogenetic network reconstruction 
contains numerous results on how, and under which conditions, a network can be reconstructed 
from triplets or suitable generalizations thereof \cite{HRS:11,JANSSON200660,vanIersel2011,doi:10.1142/S0219720012500138,5438983,Gambette2017,doi:10.1137/S0097539704446529, 10.1371/journal.pone.0106531,VANIERSEL2022106300,Semple2021,vanIersel2014}. 
However, the notion of a triplet is more subtle in the context of networks. 
To be more precise, for a DAG $G$, there are (at least) two ways to express that a triplet $ab|c$ is displayed by $G$:
\begin{enumerate}
\item[\axiom{T1}] $a,b,c \in X$ and $\lca_G(ab) \prec_G \lca_G(ac) = \lca_G(bc)$;
\item[\axiom{T2}] $a,b,c \in X$ and $G$ contains distinct vertices $u$ and $v$, as well as pairwise 
    internally vertex-disjoint directed paths $u \leadsto a$, $u \leadsto b$, $v \leadsto u$, and $v \leadsto c$.
\end{enumerate}
Note that the triplet $t = ab|c$ under \axiom{T1} naturally translates to the LCA constraint 
$R_t = \{(ab,ac), (ab,bc), (ac,bc), (bc,ac)\}$, in which case the canonical network $\cN_{R_t}$ is just the triplet $t$. 
While in rooted phylogenetic trees the conditions \axiom{T1} and \axiom{T2} are equivalent, this is not the case for general DAGs or networks (see Figure~\ref{fig:tiplesVSlca}). 
There are ways to infer triplets according to \axiom{T1} from genomic data (see e.g. ~\cite{Stadler2020,SCHALLER202163}). In this way one could, in principle, 
infer networks from genomic data by first inferring triplets according to \axiom{T1}, then creating the corresponding LCA constraints, and finally computing the canonical network in case the constraints are realizable.  It is not hard to see that triplets
satisfying \axiom{T1} also satisfy \axiom{T2}. The converse, however, 
is in general not satisfied, see Figure~\ref{fig:tiplesVSlca} for several examples.
Inferring
\axiom{T2}-triplets that do not satisfy \axiom{T1} from genomic data
will probably be more challenging. 
This is because such triplets require evidence of distinct, 
internally disjoint ancestral paths and therefore depend on explicitly 
reconstructed network structures.
\REV{Nevertheless, most algorithmic work on triples and phylogenetic networks uses the notion of \axiom{T2}-triples. Although every set $\mathscr R$ of \axiom{T2}-triples can be displayed by some network,
such networks may be uninformative, since they can display all three rooted triples on the same three leaves; see, for example, $N_3$ in Figure~\ref{fig:tiplesVSlca} and \cite{Huson:11,JANSSON200660}. It is therefore natural to seek more restricted, and hence more informative, networks displaying all \axiom{T2}-triples in $\mathscr R$.
However, imposing structural constraints often makes the corresponding decision problems NP-hard \cite{doi:10.1137/S0097539704446529,vanIersel2009}. By contrast, polynomial-time algorithms are known for constructing particular types of networks when strong restrictions are imposed on the set $\mathscr R$ of \axiom{T2}-triples; see, e.g., \cite{Aho:81,vanIersel2011,lev2-09,5438983}.

To the best of our knowledge, no concise framework has yet been established for constructing DAGs or phylogenetic networks that display a given set of \axiom{T1}-triples. We believe that the theory of realizable relations may provide a foundation for such a framework. This will be explored in future work.
}

\paragraph{Limitation and Extension of Condition \axiom{X1} and \axiom{X2}.}
The implication in Equation~\eqref{eq:LCA-implication} of
Lemma~\ref{lem:simpleLCAinfer} is only one of many that can be derived from the ``LCA-$\prec_G$ relationships''.
Nevertheless, Equation~\eqref{eq:LCA-implication} laid the foundation for the notion of 
cross-consistency, leading to the $+$-closure $R^+$ of a relation $R$. 
This, in turn, allowed the definition of conditions \axiom{X1} and \axiom{X2}, 
which characterize relations that are realizable by some DAG.
However, for practical purposes it would be interesting to understand when a relation can be realized by a DAG that lies within a fixed class.
There are examples of subsets $R \subseteq \rel_H$, with $H$ belonging to a class $\mathscr{N}$ 
of DAGs, for which $\cG_R$ does not yield a network in $\mathscr{N}$. 
This occurs mainly because $R$ contains only partial information about $\rel_H$, see Figure~\ref{fig:non-min}. 
In this example, the vertices $[ab]$ and $[bc]$ in the canonical network $\cN_R=\cG_R$
are never identified by our approach, 
since $ab \, R^+ \, bc$ and \REV{$bc \, R^+ \, ab$} does not hold. 
Moreover, since neither $ab \, R^+ \, bc$ nor $bc \, R^+ \, ab$ holds, 
 the vertices $[ab]$ and $[bc]$ are $\preceq_{\cN_R}$-incomparable.
 However, $R$ is realized by the tree $T$ shown in Figure~\ref{fig:non-min}. 
Consequently, to infer DAGs or networks belonging to a particular class $\mathscr{N}$, 
such as phylogenetic trees or other types of networks, \axiom{X1} and \axiom{X2} 
must be extended by adding properties that will ensure that a DAG within 
$\mathscr{N}$ is returned whenever one exists that realizes $R$. 
If $\mathscr{N}$ is the class of phylogenetic trees, the \textsc{Build} algorithm could be applied; 
however, our focus is on a non-recursive approach, with particular interest in the 
structural properties of relations $R$ that are realizable by a tree. 
Similar questions can be asked for the tuple $(R,S)$, where $S$ contains additional information about incomparable LCA constraints.
A further question arising in this context is that of ``encoding'' of networks, i.e., 
for which network classes $\mathscr{N}$ 
does it hold that $\rel_N = \rel_{N'}$ if and only if $N \sim N'$ for all $N, N' 
 \in \mathscr{N}$, where $\sim$ denotes an isomorphism between $N$ and ${N'}$
 that is the identity on the leaf-set?

\paragraph{Underlying Optimization Problems.}
The are various interesting optimization problems that arise in the
context of our results.
One problem concerns finding the ``simplest'' DAG or network that realizes a given relation $R$, 
a notion that can be formalized in various  ways. 
For instance, if one restricts to trees, a natural objective 
could be to find a tree with the minimum number of vertices among all trees realizing $R$, a problem that is known to be NP-hard~\cite{JLL:12}. 
In our setting, the canonical DAG $\cG_R$ is not necessarily minimum-sized, i.e., it does not always have the fewest vertices among all DAGs realizing $R$ (see for example Figures~\ref{fig:non-min} and \ref{fig:exmpl-tree-reg}). 
Determining the computational complexity of finding a DAG or network realizing $R$ with the minimum number of vertices thus remains an open problem. 
However, in DAGs, ``simplicity'' need not be necessarily measured by vertex count. 
An alternative criterion is proximity to a tree structure, for instance, minimizing the total number of hybrid vertices (those of in-degree 
greater than one)
or the maximum number of hybrids within any biconnected component\footnote{
A \emph{biconnected component} of $G$ is an inclusion-maximal subgraph $H=(V',E')$ of $G$ such that $H$ remains weakly connected after removal of any one of its vertices and its incident arcs \cite{book:91867427}.}. 
We therefore also pose the problem of determining the computational complexity of finding a DAG or network that realizes $R$ with as few hybrids as possible, either globally or per biconnected component.
Similar questions can be asked for tuples $(R,S)$.

A third optimization problem concerns non-realizable relations. 
Recall that only relations satisfying \axiom{X1} and \axiom{X2} are realizable (cf.~Theorem~\ref{thm:char}). 
Given a non-realizable relation $R$, one may ask whether it is possible to efficiently determine the largest subset $R'\subseteq R$ that is realizable, or whether this problem is NP-hard. 
If $R$ is a set of LCA constraints or represents a set of triplets, it has been shown that finding a maximum-sized subset that can be realized by a tree $T$ is NP-hard~\cite{10.1007/BFb0045080,JANSSON200150,Wu2004}. 
However, our framework allows greater flexibility: we can consider sets $R$ containing LCA constraints (respectively, LCA constraints of triplets according to \axiom{T1}), 
and ask when they can be realized by some DAG or network without restricting ourselves to trees. 
As shown in Proposition~\ref{prop:CanonicalG-DAG-outsourced} and Example~\ref{exmpl:RnotR+is}, there exist non-realizable relations $R$ for which $R^+$ is realizable, 
in particular, those satisfying \axiom{X1} but not \axiom{X2}. 
This raises the question of whether the canonical DAG $\cG_R$ that realizes $R^+$ 
could serve as a scaffold for determining a largest realizable subset, or for identifying minimal conflicts within relations satisfying \axiom{X1}. 
In other words, can $\cG_R$ be used to define a measure of how far a relation is from being realizable
which can be exploited in practice? 
Or is $\cG_R$ possibly already a DAG that realizes as many constraints in $R$ as possible?
Similar questions can be asked for tuples $(R,S)$.

\begin{figure}
    \centering
    \includegraphics[width=0.8\textwidth]{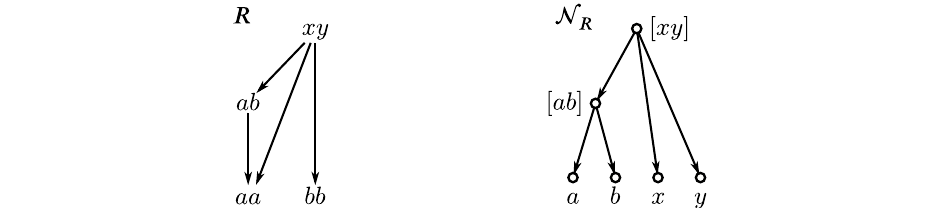}
    \caption{On the left we give a graphical representation of a relation relation $R = \{(aa,xy), (bb,xy), (ab,xy), (aa,ab)\}$
    on $\pairs(X)$ with $X =\{a,b,x,y\}$. Here we draw an arc $p\to q$ precisely if $q\ R\ p$.
    On the right the canonical network $\cN_R$ is shown. 
    Both subsets $R'=\{(ab,xy)\}$ and $R''=\{(aa,xy), (bb,xy),  (aa,ab)\}$ are inclusion-minimal
    w.r.t.\ the property that $\cl(R')=\cl(R'')=\cl(R)$. However, $|R'|\neq |R''|$.
    }
    \label{fig:no-matroid}
\end{figure}

\paragraph{Matroids.}

 The study of matroids in phylogenetics has become an active research area in recent years \cite{AK:06,DHS:14,Francisco2009,HOLLERING2021142,MaayanLevy24}. 
 A classical question in this context is to find ``minimum closure-representatives'', i.e., minimum-sized subsets $R' \subseteq R$ 
 where $R$ encodes partial information about phylogenetic networks or trees and $R'$ satisfies $\mathrm{CL}(R') = \mathrm{CL}(R)$ 
 for a well-defined closure operator $\mathrm{CL}$ \cite{SH:18,MaayanLevy24}. 
\REV{For example, let $\mathscr R$ be a set of triples that is compatible with a tree $T$, i.e., all triples in $\mathscr R$ are displayed by $T$ according to \axiom{T1}. Let $\mathfrak T_{\mathscr R}$ denote the set of all trees compatible with $\mathscr R$, and let $\mathscr R(T)$ be the set of triples displayed by $T$
according to \axiom{T1}. Then the intersection
\[
\mathrm{CL}(\mathscr R)=\bigcap_{T\in \mathfrak T_{\mathscr R}} \mathscr R(T)
\]
contains precisely those triples that are displayed by every tree compatible with $\mathscr R$. In this setting, the minimum-size subsets $\mathscr R'$ with $\mathrm{CL}(\mathscr R')=\mathrm{CL}(\mathscr R)$ form the bases of a matroid and can be found by a simple greedy algorithm~\cite{SH:18}. As shown by M.\ Levy, these type of matroids are graphic, i.e., 
they are isomorphic to the cycle matroid of some undirected graph \cite{MaayanLevy24}.}
 In our setting of relations, this corresponds to finding minimum-sized subsets $R' \subseteq R$ such that $\cl(R') = \cl(R)$. 
 However, as illustrated in Figure~\ref{fig:no-matroid}, it can happen that $|R'| \neq |R''|$ for different minimum-sized subsets, 
 showing that, in general, no matroid structure underlies LCA constraints and their closure.
The latter raises the question of how difficult it is to find a minimum-sized subset $R' \subseteq R$ 
for a realizable relation $R$ such that $\cl(R') = \cl(R)$. 
Moreover, one may restrict attention to relations $R$ that can be realized by particular classes of networks, 
for example, phylogenetic trees. 
This leads to the question of whether there exist network classes $\mathscr{N}$ for which the tuple 
\[
(R, \mathcal{F}_R), \quad \text{where } \mathcal{F}_R = \{ R'' \subseteq R' \mid R' \subseteq R \text{ is inclusion-minimal with } \cl(R') = \cl(R) \},
\]
forms a matroid, with $R$ a relation realized by some network $N \in \mathscr{N}$.
Moreover, if we restrict $R$ to relations that contain ``full information'' about the underlying network, 
i.e., $R = \rel_N$ for some network $N \in \mathscr{N}$, 
does this make it easier to determine minimum-sized subsets $R' \subseteq R$ such that $\cl(R') = \cl(\rel_N) = \rel_N$?
Similar questions can be asked for tuples $(R,S)$.

\section*{Acknowledgment}

We thank the organizers of the Bielefeld meeting ``New Directions in Experimental Mathematics'' 
(June 30 -- July 1, 2025), where MH proposed the open problem on LCA constraints for DAGs and networks, 
and the BIRS Workshop ``Novel Mathematical Paradigm for Phylogenomics'' 
(August 24--29, 2025), where further discussions helped refine the ideas presented in this work.

We also thank  Nadia Tahiri, Leo van Iersel,  Peter F.\ Stadler and  Lars Berling for fruitful
discussions on the topic. 
Finally, we also thank two anonymous referees for their helpful comments.

\printcredits

\section*{Declaration of competing interest}
The authors declare that they have no known competing financial interests or personal relationships that could have appeared to influence the work reported in this paper.

\bibliographystyle{cas-model2-names}

\bibliography{common}


\end{document}